\documentclass[11pt]{article}
\usepackage{mathrsfs}
\usepackage{amscd}
\usepackage{amsmath,amsfonts,amssymb,amscd}
\usepackage{indentfirst,graphics,epsfig,psfrag}
\input{epsf}
\usepackage{ifpdf}
\usepackage{enumerate}
\usepackage{appendix}
\usepackage{enumerate}
\usepackage{color}
\usepackage{lineno}

\usepackage{lineno}
\makeatletter \@addtoreset{figure}{section} \makeatother
\makeatletter
\long\def\@makecaption#1#2{%
   \vskip 10\p@
   \setbox\@tempboxa\hbox{{#1}\ \ #2}%
   \ifdim \wd\@tempboxa >\hsize

       {#1}\ \ #2\par
   \else
       \hbox to\hsize{\hfil\box\@tempboxa\hfil}%
   \fi}
\makeatother

\newtheorem{thm}{Theorem}[section]
\newtheorem{cor}[thm]{Corollary}
\newtheorem{lem}{Lemma}[section]

\newtheorem{obs}[thm]{Observation}
\newtheorem{pro}[thm]{Proposition}

\newcommand{\qed}{{\hfill\rule{3pt}{7pt}}}

\def\qed{\hfill \rule{4pt}{7pt}}

\begin{document}
\title{{\bf Line $k$-Arboricity in Product Networks}
\footnote{Supported by the National Science Foundation of China
(Nos. 11551001, 11161037, 11461054) and the Science Found of Qinghai
Province (No. 2014-ZJ-907).}}
\author{
\small Yaping Mao$^{1,2}$, \ \ Zhiwei Guo$^1$,\footnote{Corresponding author} \ \ Nan Jia$^1$, \ \ He Li$^1$\\[0.2cm]
\small $^1$Department of Mathematics, Qinghai Normal University,\\
\small Xining, Qinghai 810008, China\\[0.2cm]
\small $^2$Key Laboratory of IOT of Qinghai Province,\\
\small Xining, Qinghai 810008, China\\[0.2cm]
\small E-mails: maoyaping@ymail.com; guozhiweic@yahoo.com\\
\small ~~~~~~~~~~~jianan66ss@yahoo.com; lihe0520@yahoo.com\\
\small }
\date{}
\maketitle
\begin{abstract}
A \emph{linear $k$-forest} is a forest whose components are paths of
length at most $k$. The \emph{linear $k$-arboricity} of a graph $G$,
denoted by ${\rm la}_k(G)$, is the least number of linear
$k$-forests needed to decompose $G$. Recently, Zuo, He and Xue
studied the exact values of the linear $(n-1)$-arboricity of
Cartesian products of various combinations of complete graphs,
cycles, complete multipartite graphs. In this paper, for general $k$
we show that $\max\{{\rm la}_{k}(G),{\rm la}_{\ell}(H)\}\leq {\rm
la}_{\max\{k,\ell\}}(G\Box H)\leq {\rm la}_{k}(G)+{\rm
la}_{\ell}(H)$ for any two graphs $G$ and $H$. Denote by $G\circ H$,
$G\times H$ and $G\boxtimes H$ the lexicographic product, direct
product and strong product of two graphs $G$ and $H$, respectively.
We also derive upper and lower bounds of ${\rm la}_{k}(G\circ H)$,
${\rm la}_{k}(G\times H)$ and ${\rm la}_{k}(G\boxtimes H)$ in this
paper. The linear $k$-arboricity of a $2$-dimensional grid graph, a
$r$-dimensional mesh, a $r$-dimensional torus, a $r$-dimensional
generalized hypercube and a $2$-dimensional hyper Petersen network are also studied.\\[2mm]
{\bf Keywords:} Linear $k$-forest, linear $k$-arboricity, Cartesian
product, complete product, lexicographical product, strong product, direct product.\\[2mm]
{\bf AMS subject classification 2010:} 05C15, 05C76, 05C78.
\end{abstract}

\section{Introduction}

All graphs considered in this paper are undirected, finite and
simple. We refer to the book \cite{Bondy} for graph theoretical
notation and terminology not described here. Let $N$ be the set of
natural numbers and let $[a,b]$ be the set $\{n\in N\mid a\leq n\leq
b\}$. A \emph{decomposition} of a graph is a list of subgraphs such
that each edge appears in exactly one subgraph in the list. If a
graph $G$ has a decomposition $G_1,G_2,\ldots,G_t$, then we say that
\emph{$G_1,G_2,\ldots,G_t$ decompose $G$} or \emph{$G$ can be
decomposed into $G_1,G_2,\ldots,G_t$}. Furthermore, a \emph{linear
$k$-forest} is a forest whose components are paths of length at most
$k$. The \emph{linear $k$-arboricity} of a graph $G$, denoted by
${\rm la}_k(G)$, is the least number of linear $k$-forests needed to
decompose $G$.

The notion of linear $k$-arboricity of a graph was first introduced
by Habib and Peroche \cite{Habib}, which is a natural generalization
of edge-coloring. Clearly, a linear $1$-forest is induced by a
matching, and ${\rm la}_1(G)$ is the chromatic index $\chi'(G)$ of a
graph $G$. Moreover, the linear $k$-arboricity ${\rm la}_k(G)$ is
also a refinement of the ordinary linear arboricity ${\rm la}(G)$
(or ${\rm la}_{\infty}(G)$ which is the case when every component of
each forest is a path with no length constraint. By the way, the
notion of linear arboricity was introduced earlier by Harary in
\cite{Harary}. For more details on linear $k$-arboricity, we refer
to \cite{AldredW, Alon, AlonTW, Chang, ChangCFH, ChenH}

In graph theory, Cartesian product, strong product, lexicographical
product and direct product are four of main products, each with its
own set of applications and theoretical interpretations. Product
networks were proposed based upon the idea of using the cross
product as a tool for ``combining'' two known graphs with
established properties to obtain a new one that inherits properties
from both \cite{DayA}. Recently, there has been an increasing
interest in a class of interconnection networks called Cartesian
product networks; see \cite{Bao, DayA, Ku}.

The join, Cartesian, lexicographical, strong and direct products are
defined as follows.
\begin{itemize}
\item The {\it join} or {\it complete product} $G\vee H$ of two
disjoint graphs $G$ and $H$, is the graph with vertex set $V(G)\cup
V(H)$ and edge set $E(G)\cup E(H)\cup \{uv\,|\, u\in V(G), v\in
V(H)\}$.

\item The {\it Cartesian product} $G\Box H$ of two graphs $G$ and $H$, is
the graph with vertex set $V(G)\times V(H)$, in which two vertices
$(u,v)$ and $(u',v')$ are adjacent if and only if $u=u'$ and
$(v,v')\in E(H)$, or $v=v'$ and $(u,u')\in E(G)$.

\item The {\it lexicographic product} $G\circ H$ of graphs $G$ and $H$ has
the vertex set $V(G\circ H)=V(G)\times V(H)$, and two vertices
$(u,v),(u',v')$ are adjacent if $uu'\in E(G)$, or if $u=u'$ and
$vv'\in E(H)$.

\item The \emph{strong product} $G\boxtimes H$ of graphs $G$ and
$H$ has the vertex set $V(G)\times V(H)$. Two vertices $(u,v)$ and
$(u',v')$ are adjacent whenever $uu'\in E(G)$ and $v=v'$, or $u=u'$
and $vv' \in E(H)$, or $uu'\in E(G)$ and $vv'\in E(H)$.

\item The \emph{direct product} $G\times H$ of graphs $G$ and
$H$ has the vertex set $V(G)\times V(H)$. Two vertices $(u,v)$ and
$(u',v')$ are adjacent if the projections on both coordinates are
adjacent, i.e., $uu'\in E(G)$ and $vv'\in E(H)$.
\end{itemize}

From the definition of linear $k$-arboricity and the structure of
graph product, the following result is immediate.
\begin{obs}\label{obs1-1}
Let $G$ and $H$ be two graphs. Then
$$
{\rm la}_{k}(G\boxtimes H)\leq {\rm la}_{k}(G\Box H)+{\rm
la}_{k}(G\times H).
$$
\end{obs}

Xue and Zuo \cite{XueZ} investigated the linear $(n-1)$-arboricity
of complete multipartite graphs. Recently, Zuo, He and Xue
\cite{ZuoHX} studied the exact values of the linear
$(n-1)$-arboricity of Cartesian products of various combinations of
complete graphs, cycles, complete multipartite graphs.

In this paper, we consider four standard products: the
lexicographic, the strong, the Cartesian and the direct with respect
to the linear $k$-arboricity. Every of these four products will be
treated in one of the forthcoming subsections in Section 2. In
Section $3$, we demonstrate the usefulness of the proposed
constructions by applying them to some instances of product
networks.

\section{Results for general graphs}

As usual, the \emph{union} of two graphs $G$ and $H$ is the graph,
denoted by $G\cup H$, with vertex set $V(G)\cup V(H)$ and edge set
$E(G)\cup E(H)$. The disjoint union of $k$ copies of the same graph
$G$ is denoted by $kG$. The \emph{join} $G\vee H$ of two disjoint
graphs $G$ and $H$ is obtained from $G\cup H$ by joining each vertex
of $G$ to every vertex of $H$. In the sequel, let $K_{s,t}$,
$C_{n}$, $K_{n}$ and $P_n$ denote the complete bipartite graph of
order $s+t$ with part sizes $s$ and $t$, cycle of order $n$,
complete graph of order $n$, and path of order $n$, respectively.

In the sequel, let $G$ and $H$ be two connected graphs with
$V(G)=\{u_1,u_2,\ldots,u_{n}\}$ and $V(H)=\{v_1,v_2,\ldots,v_{m}\}$,
respectively. Then $V(G* H)=\{(u_i,v_j)\,|\,1\leq i\leq n, \ 1\leq
j\leq m\}$, where $*$ denotes a kind of graph product operations.
For $v\in V(H)$, we use $G(v)$ to denote the subgraph of $G* H$
induced by the vertex set $\{(u_i,v)\,|\,1\leq i\leq n\}$.
Similarly, for $u\in V(G)$, we use $H(u)$ to denote the subgraph of
$G* H$ induced by the vertex set $\{(u,v_j)\,|\,1\leq j\leq m\}$.

The following observations are immediate.
\begin{obs}\label{obs2-1}
Let $H$ be a subgraph of $G$. If $\ell\geq k$, then
$$
{\rm la}_{\ell}(H)\leq {\rm la}_{k}(G).
$$
\end{obs}
\begin{obs}{\upshape \cite{Yeh}}\label{obs2-2}
If a graph $G$ is the edge-disjoint union of two subgraphs $G_1$ and
$G_2$, then
$$
{\rm la}_{k}(G)\leq {\rm la}_{k}(G_1)+{\rm la}_{k}(G_2).
$$
\end{obs}
\begin{obs}{\upshape \cite{Yeh}}\label{obs2-3}
If a graph $G$ is the disjoint union of two subgraphs $G_1$ and
$G_2$, then
$$
{\rm la}_{k}(G)=\max\{{\rm la}_{k}(G_1),{\rm la}_{k}(G_2)\}.
$$
\end{obs}
\begin{obs}\label{obs2-4}
If $G$ is not a forest, then ${\rm la}_{k}(G)\geq 2$ for $k\geq 1$.
\end{obs}

\begin{lem}{\upshape \cite{ChenH}}\label{lem2-1}
For a graph $G$ of order $n$,
$$
\Delta(G)+1\geq \chi'(G)={\rm la}_{1}(G)\geq {\rm la}_{2}(G)\geq
\ldots \geq {\rm la}_{n-1}(G)={\rm la}(G),
$$
where $\chi'(G)$ denotes the edge chromatic number of $G$.
\end{lem}

\begin{lem}{\upshape \cite{ChenH}}\label{lem2-2}
For a graph $G$,
$$
{\rm la}_{k}(G)\geq
\max\left\{\left\lceil\frac{\Delta(G)}{2}\right\rceil, \left\lceil
\frac{|E(G)|}{\lfloor\frac{k|V(G)|}{k+1}\rfloor}\right\rceil\right\}.
$$
\end{lem}

\subsection{For Cartesian product}

We first give the bounds for general graphs.
\begin{thm}\label{th2-5}
Let $G$ and $H$ be two graphs. Then
$$
\max\{{\rm la}_{k}(G),{\rm la}_{\ell}(H)\}\leq {\rm
la}_{\max\{k,\ell\}}(G\Box H)\leq {\rm la}_{k}(G)+{\rm
la}_{\ell}(H).
$$
Moreover, the upper bound is sharp.
\end{thm}
\begin{pf}
Set $V(G)=\{u_1,u_2,\ldots,u_{n}\}$ and
$V(H)=\{v_1,v_2,\ldots,v_{m}\}$. Let ${\rm la}_{k}(G)=p$ and ${\rm
la}_{k}(H)=q$. Since ${\rm la}_{k}(G)=p$, it follows that there are
$p$ linear $k$-forests in $G$. Then, since ${\rm la}_{k}(G(v_i))=p$,
it follows that there are $p$ linear $k$-forests in $G(v_i)$, say
$F_{i,1},F_{i,2},\ldots,F_{i,p}$. Similarly, since ${\rm
la}_{k}(H)=q$, it follows that there are $q$ linear $k$-forests in
$H$. Then, since ${\rm la}_{k}(H(u_i))=q$, it follows that there are
$q$ linear $k$-forests in $H(u_i)$, say
$F_{i,1}',F_{i,2}',\ldots,F_{i,q}'$. Set
$F_j=\bigcup_{i=1}^mF_{i,j}$ where $1\leq j\leq p$, and
$F_j'=\bigcup_{i=1}^nF_{i,j}'$ where $1\leq j\leq q$. Clearly,
$F_1,F_2,\ldots,F_p,F_1',F_2',\ldots,F_q'$ are $p+q$ linear
$k$-forest in $G\Box H$. So ${\rm la}_{\max\{k,\ell\}}(G\Box H)\leq
{\rm la}_{k}(G)+{\rm la}_{\ell}(H)$. From Observation \ref{obs2-1},
we have ${\rm la}_{\max\{k,\ell\}}(G\Box H)\geq \max\{{\rm
la}_{k}(G),{\rm la}_{\ell}(H)\}$.\qed
\end{pf}\\

The following corollary is a generalization of the above result.
\begin{cor}\label{cor2-6}
Let $G_1,G_2,\ldots,G_r$ be graphs. Then
\begin{eqnarray*}
\max\{{\rm la}_{k_1}(G_1),{\rm la}_{k_2}(G_2),\ldots,{\rm
la}_{k_r}(G_r)\}&\leq& {\rm la}_{\max\{k_1,k_2,\ldots,k_r\}}(G_1\Box
G_2\Box \ldots \Box G_r)\\
&\leq&{\rm la}_{k_1}(G_1)+{\rm la}_{k_2}(G_2)+\ldots+{\rm
la}_{k_r}(G_r).
\end{eqnarray*}

Moreover, the bounds are sharp.
\end{cor}

\subsection{For complete product}

The following results were obtained by Dirac \cite{Dirac}; see
Laskar and Auerbach \cite{LA}.
\begin{pro}{\upshape\cite{Dirac, LA}} \label{pro2-7}
$(1)$ For all even $r\geq 2$, $K_{r,r}$ is the union of its
$\frac{1}{2}r$ Hamiltonian cycles.

$(2)$ For all odd $r\geq 3$, $K_{r,r}$ is the union of its
$\frac{1}{2}r$ Hamiltonian cycles and one perfect matching.
\end{pro}

For complete product, we have the following.
\begin{thm}\label{th2-8}
Let $G$ and $H$ be two graphs. Then
$$
\max\left\{\left\lceil
\frac{\Delta(G)+|V(H)|}{2}\right\rceil,\left\lceil
\frac{\Delta(H)+|V(G)|}{2}\right\rceil\right\}\leq {\rm
la}_{k}(G\vee H)\leq {\rm la}_{k}(G)+{\rm la}_{k}(H)+|V(H)|.
$$
\end{thm}
\begin{pf}
Set $|V(G)|=n$ and $|V(H)|=m$. Without loss of generality, let
$n\leq m$. Let $G'=G\cup (m-n)K_1$. Then $|V(G')|=m$ and $G\vee H$
is a subgraph of $G'\vee H$. Since ${\rm la}_{k}(G)=p$, it follows
that there are $p$ linear $k$-forests in $G$, say
$F_1,F_2,\ldots,F_p$. From the structure of $G'$,
$F_1,F_2,\ldots,F_p$ are linear $k$-forests in $G'$, and hence ${\rm
la}_{k}(G')\leq p$. From Observation \ref{obs2-1}, we have $p={\rm
la}_{k}(G)\leq {\rm la}_{k}(G')\leq p$, and hence ${\rm
la}_{k}(G')=p$. Since ${\rm la}_{k}(H)=q$, it follows that there are
$q$ linear $k$-forests in $H$, say $F_1',F_2',\ldots,F_q'$. Note
that the subgraph induced by all the vertices of $G'\vee H$ is a
complete bipartite graph $K_{m,m}$. From Proposition \ref{pro2-7},
$K_{m,m}$ can be decomposed into $m$ perfect matchings, say
$M_1,M_2,\ldots,M_m$. These perfect matchings are $m$ linear
$1$-forests in $G'\vee H$. Observe that
$$
E(G'\vee H)=\left(\bigcup_{i=1}^mM_i\right)\cup
\left(\bigcup_{i=1}^pE(F_i)\right)\cup
\left(\bigcup_{i=1}^qE(F_i')\right).
$$
Then $F_1,F_2,\ldots,F_p,F_1',F_2',\ldots,F_q',M_1,M_2,\ldots,M_m$
form $p+q+m$ linear $k$-forests in $G'\vee H$. So ${\rm
la}_{k}(G\vee H)\leq {\rm la}_{k}(G'\vee H)\leq {\rm la}_{k}(G)+{\rm
la}_{k}(H)+|V(H)|$. Note that $\Delta(G\vee H)\geq
\max\{\Delta(G)+|V(H)|,\Delta(H)+|V(G)|\}$. From Lemma \ref{lem2-2},
we have
$$
{\rm la}_{k}(G\vee H)\geq \max\left\{\left\lceil
\frac{\Delta(G)+|V(H)|}{2}\right\rceil,\left\lceil
\frac{\Delta(H)+|V(G)|}{2}\right\rceil\right\},
$$
as desired. \qed
\end{pf}

\subsection{For lexicographical product}

From the definition, the lexicographic product graph $G\circ H$ is a
graph obtained by replacing each vertex of $G$ by a copy of $H$ and
replacing each edge of $G$ by a complete bipartite graph $K_{m,m}$.
For an edge $e=u_iu_j\in E(G) \ (1\leq i,j\leq n)$, the induced
subgraph obtained from the edges between the vertex set
$V(H(u_i))=\{(u_i,v_1),(u_i,v_2),\cdots,(u_i,v_{m})\}$ and the
vertex set $V(H(u_j))=\{(u_j,v_1),(u_j,v_2),\cdots,$ $(u_j,v_{m})\}$
in $G\circ H$ is a complete equipartition bipartite graph of order
$2m$, denoted by $K_{e}$ or $K_{u_i,u_j}$.

From Proposition \ref{pro2-7}, $K_e$ can be decomposed into $m$
perfect matching, denoted by $M^e_1,M^e_2,\ldots,M^e_{m}$. We now in
a position to give the result for lexicographical product.
\begin{thm}\label{th2-9}
Let $G$ and $H$ be two graphs. Then
$$
\left\lceil \frac{\Delta(H)+|V(H)|\Delta(G)}{2}\right\rceil\leq {\rm
la}_{\max\{k,\ell\}}(G\circ H)\leq {\rm la}_{k}(G)|V(H)|+{\rm
la}_{\ell}(H).
$$

Moreover, the bounds are sharp.
\end{thm}
\begin{pf}
Set $V(G)=\{u_1,u_2,\ldots,u_{n}\}$ and
$V(H)=\{v_1,v_2,\ldots,v_{m}\}$. For $k\leq \ell$, we need to show
that ${\rm la}_{\ell}(G\circ H)\leq {\rm la}_{k}(G)|V(H)|+{\rm
la}_{\ell}(H)$. Let ${\rm la}_{k}(G)=p$ and ${\rm la}_{\ell}(H)=q$.
Since ${\rm la}_{\ell}(H)=q$, it follows that there are $q$ linear
$\ell$-forest in $H$. Then, since ${\rm la}_{\ell}(H(u_i))=q$, it
follows that there are $q$ linear $\ell$-forests in $H(u_i)$, say
$F_{i,1}',F_{i,2}',\ldots,F_{i,q}'$. Set
$F_j'=\bigcup_{i=1}^qF_{i,j}'$ where $1\leq j\leq q$.

Since ${\rm la}_{k}(G)=p$, it follows that there are $p$ linear
$k$-forests in $G$, say $F_{1},F_{2},\ldots,F_{p}$. For each linear
$k$-forest $F_i \ (1\leq i\leq p)$ in $G$, we define a subgraph
$\mathcal{F}_i$ of $G\circ H$ corresponding to $F_i$ as follows:
$V(\mathcal{F}_i)=V(F_i\circ H)$ and
$E(\mathcal{F}_i)=\{(u_p,v_s)(u_q,v_t)\,|\,u_pu_q\in E(F_i), \
u_p,u_q\in V(G), \ v_s,v_t\in V(H)\}$. We call $\mathcal{F}_i$ a
\emph{blow-up linear $k$-forest corresponding to $T_i$ in $G$}; see
Figure $1$ for an example.
\begin{figure}[h,t,b,p]
\begin{center}
\scalebox{0.6}[0.6]{\includegraphics{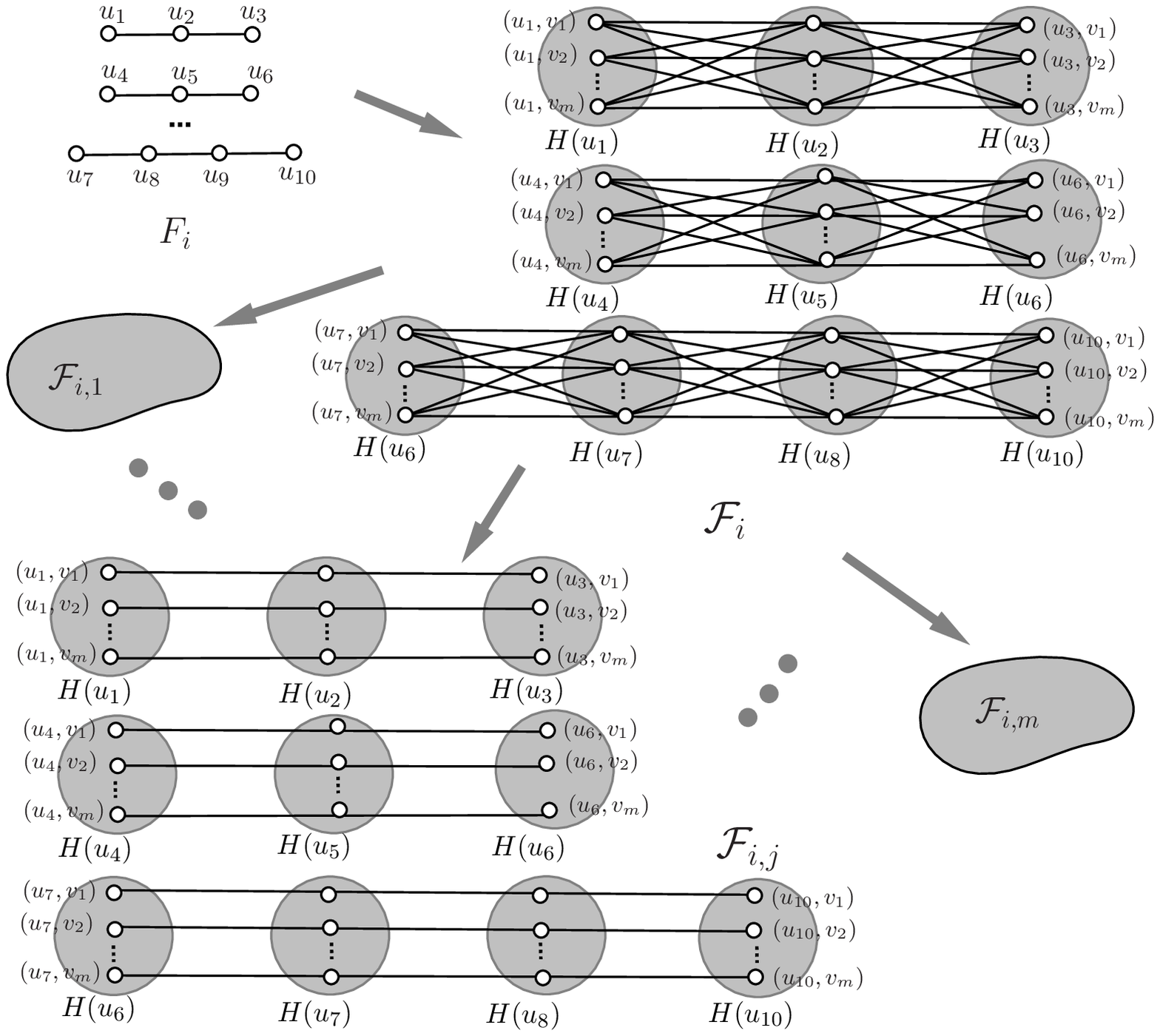}}\\
Figure 1: The blow-up linear $k$-forest $\mathcal{F}_i$ and parallel
linear $k$-forest $\mathcal{F}_{i,j}$ in $G\circ H$ corresponding to
$F_i$ in $G$.
\end{center}
\end{figure}
For each $i \ (1\leq i\leq p)$ and each $j\ (1\leq j\leq m)$, we
define another subgraph $\mathcal{F}_{i,j}$ of $G\circ H$
corresponding to $T_i$ in $G$ as follows:
$V(\mathcal{F}_{i,j})=V(F_i\circ H)$ and
$E(\mathcal{F}_{i,j})=\bigcup_{e\in E(F_i)}M^e_{i,j}$, where
$M^e_{i,j}$ is a matching of $K_e$. We call $\mathcal{F}_{i,j}$ a
\emph{parallel linear $k$-forest of $G\circ H$ corresponding to the
forest $F_i$ in $G$}; see Figure $1$ for an example. Note that all
the parallel linear $k$-forests in $\{\mathcal{F}_{i,j}\,|\,1\leq
i\leq p,\,1\leq j\leq m\}$ are $pm$ linear $k$-forests of $G\circ
H$.

Since $k\leq \ell$, it follows that the forests
$F_{1}',F_{2}',\ldots,F_{q}'$ and all the forests in
$\{\mathcal{F}_{i,j}\,|\,1\leq i\leq p,\,1\leq j\leq m\}$ form
$pm+q$ linear $\ell$-forests of $G\circ H$. Observe that each edge
of $G\circ H$ belongs to one of the above linear $\ell$-forests. So
${\rm la}_{\ell}(G\circ H)\leq {\rm la}_{k}(G)|V(H)|+{\rm
la}_{\ell}(H)$.

For $\ell\leq k$, we can prove that ${\rm la}_{k}(G\circ H)\leq {\rm
la}_{k}(G)|V(H)|+{\rm la}_{\ell}(H)$ similarly. We conclude that
${\rm la}_{\max\{k,\ell\}}(G\circ H)\leq {\rm la}_{k}(G)|V(H)|+{\rm
la}_{\ell}(H)$.

Note that $\Delta(G\circ H)\geq \Delta(H)+|V(H)|\Delta(G)$. From
Lemma \ref{lem2-2}, we have
$$
{\rm la}_{\max\{k,\ell\}}(G\circ H)\geq \left\lceil
\frac{\Delta(H)+|V(H)|\Delta(G)}{2}\right\rceil,
$$
as desired.\qed
\end{pf}

\subsection{For direct product}

For direct product, we have the following.
\begin{thm}\label{th2-10}
Let $G$ and $H$ be two graphs. Then
$$
\left\lceil \frac{\Delta(G)\Delta(H)}{2}\right\rceil\leq {\rm
la}_{\max\{k,\ell\}}(G\times H)\leq 2{\rm la}_{k}(G){\rm
la}_{\ell}(H).
$$
Moreover, the bounds are sharp.
\end{thm}
\begin{pf}
From Lemma \ref{lem2-2}, we have
$$
{\rm la}_{\max\{k,\ell\}}(G\times H)\geq \left\lceil
\frac{\Delta(G\times H)}{2}\right\rceil\geq \left\lceil
\frac{\Delta(G)\Delta(H)}{2}\right\rceil.
$$
It suffices to show that ${\rm la}_{\max\{k,\ell\}}(G\times H)\leq
2{\rm la}_{k}(G){\rm la}_{\ell}(H)$.

We now give the proof of this theorem, with a running example
(corresponding to Figure 2). From the symmetry of direct product, we
can assume $k\leq \ell$. We only need to show ${\rm
la}_{\ell}(G\times H)\leq 2{\rm la}_{k}(G){\rm la}_{\ell}(H)$. Set
$V(G)=\{u_1,u_2,\ldots,u_{n}\}$ and $V(H)=\{v_1,v_2,\ldots,v_{m}\}$.
Let ${\rm la}_{k}(G)=p$ and ${\rm la}_{\ell}(H)=q$. Since ${\rm
la}_{k}(G)=p$, it follows that there are $p$ linear $k$-forests in
$G$, say $F_1,F_2,\ldots,F_p$. For each $F_i \ (1\leq i\leq p)$, we
assume that $P_{i,1},P_{i,2},\ldots,P_{i,x}$ are all the paths in
$F_i$. Then $F_i=\bigcup_{j=1}^xP_{i,j}$ where $1\leq i\leq p$. Take
for example, let $P_{i,1}=P_2$, $P_{i,2}=P_3$, $P_{i,3}=P_4$; see
Figure 2 $(a)$. Then $F_i=P_{i,1}\cup P_{i,2}\cup P_{i,3}$.
\begin{figure}[h,t,b,p]
\begin{center}
\scalebox{0.6}[0.6]{\includegraphics{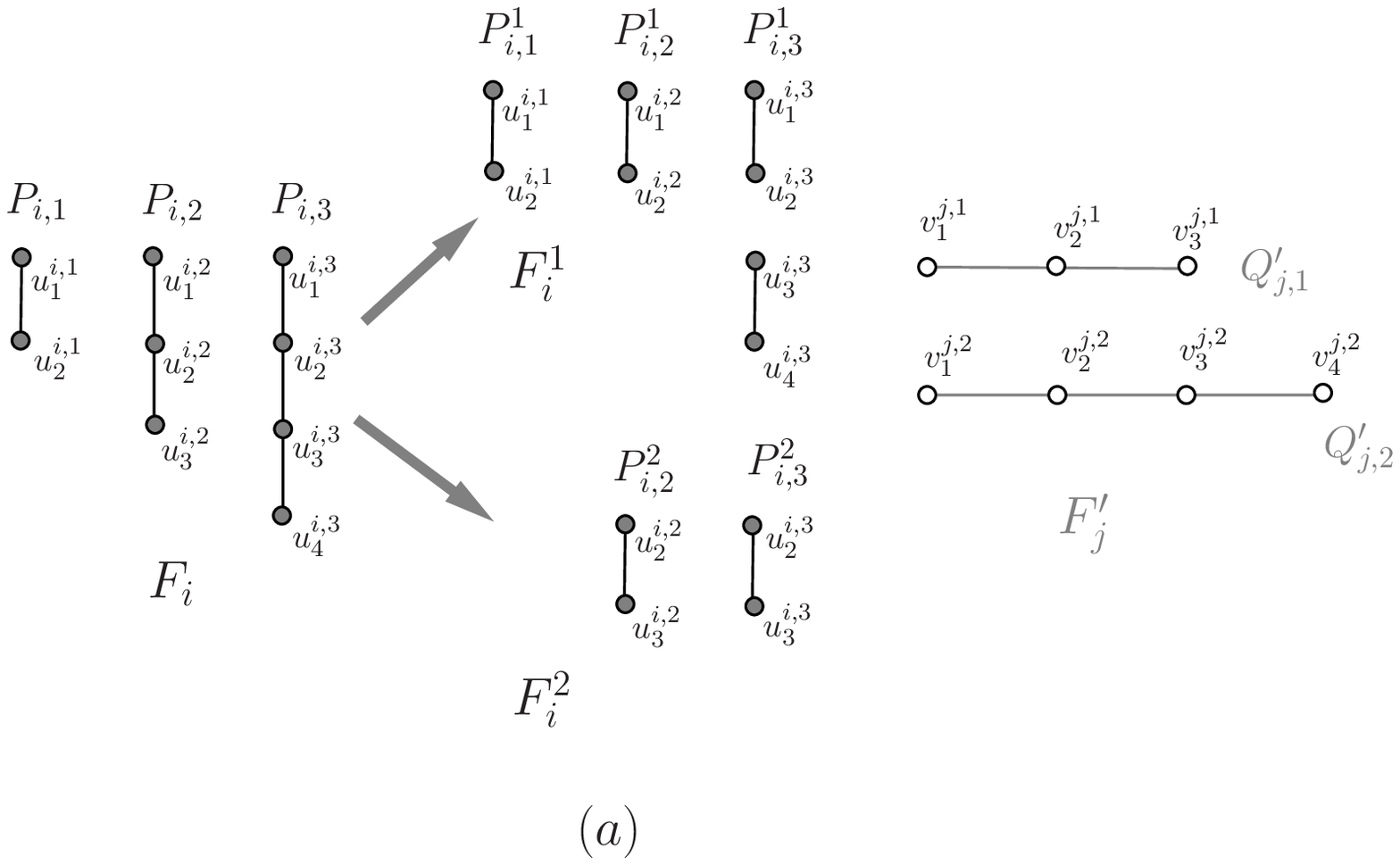}}\\[0.5cm]
\scalebox{0.6}[0.6]{\includegraphics{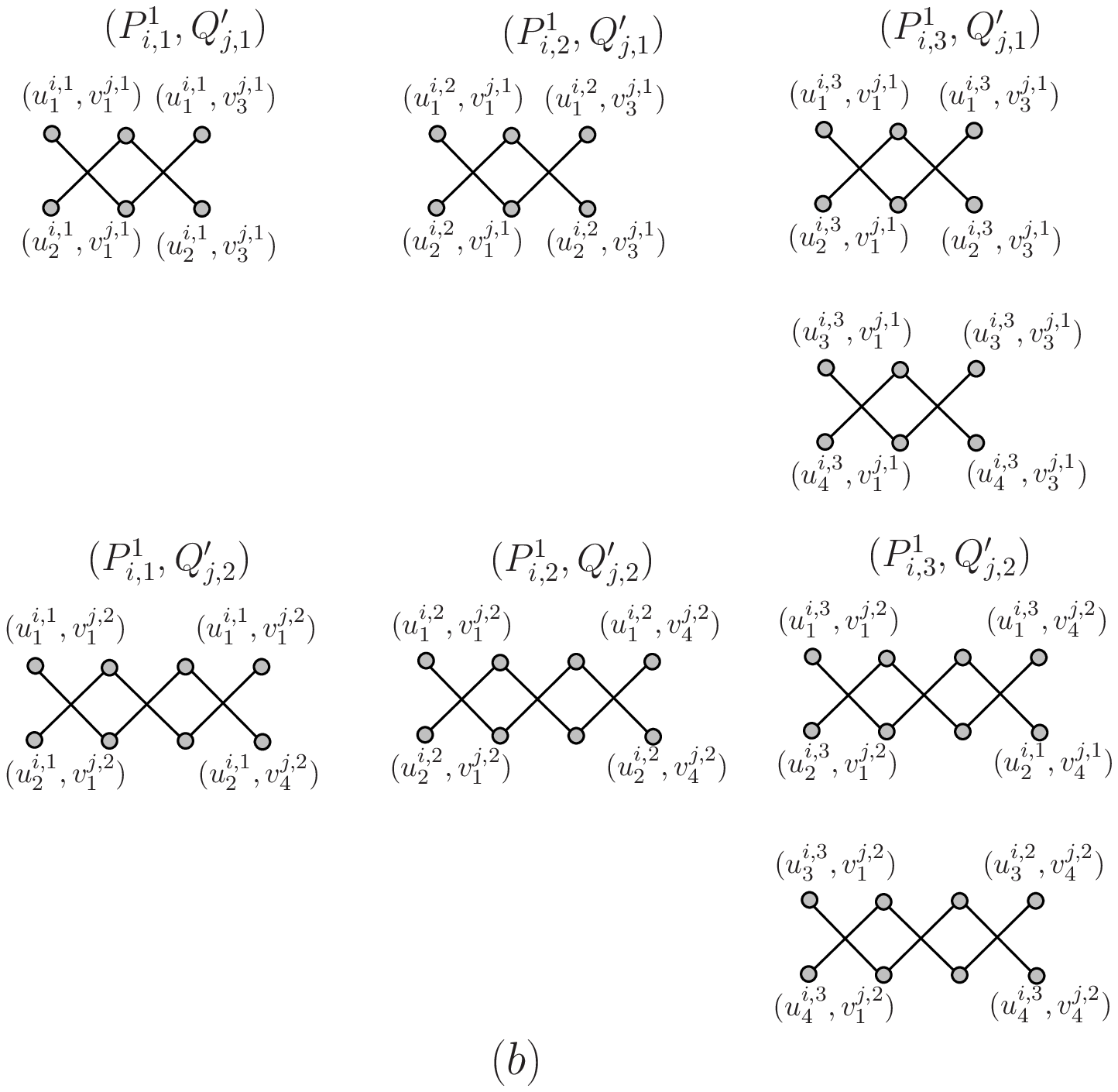}}\\[0.5cm]
\scalebox{0.6}[0.6]{\includegraphics{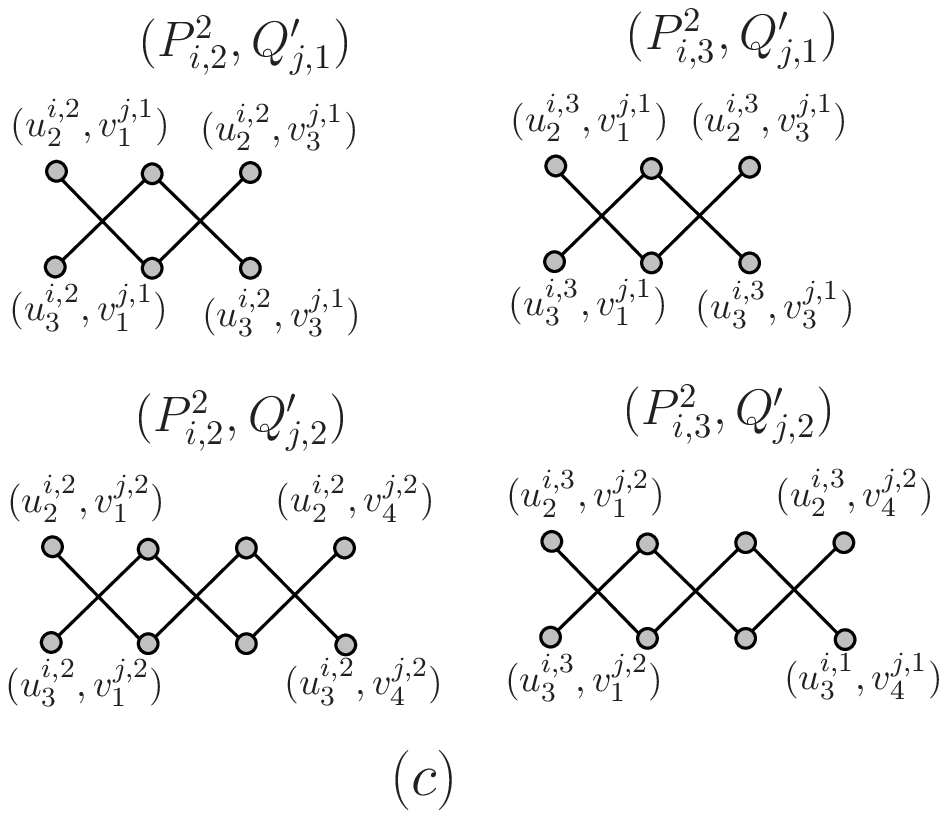}}\\
Figure 2: The running example for Theorem \ref{th2-10}.
\end{center}
\end{figure}

For each $P_{i,j}$, we let $P_{i,j}=u_1^{i,j}u_2^{i,j}\ldots
u_a^{i,j}$, where $1\leq j\leq x$. We first define two subgraphs
$P_{i,j}^1,P_{i,j}^2$ induced by the edges in
\begin{eqnarray*}
&&E(P_{i,j}^1)=\{u_{2r-1}^{i,j}u_{2r}^{i,j}\,|\,1\leq r\leq \lfloor
a/2\rfloor\},\\
&&E(P_{i,j}^2)=\{u_{2r}^{i,j}u_{2r+1}^{i,j}\,|\,1\leq r\leq \lfloor
a/2\rfloor\},
\end{eqnarray*}
respectively. Next, we set $F_i^1=\bigcup_{j=1}^xP_{i,j}^1$ and
$F_i^2=\bigcup_{j=1}^xP_{i,j}^2$. For the above example, we have
$P_{i,1}=u_1^{i,1}u_2^{i,1}$, $P_{i,2}=u_1^{i,2}u_2^{i,2}u_3^{i,2}$
and $P_{i,3}=u_1^{i,3}u_2^{i,3}u_3^{i,3}u_4^{i,3}$. Then $P_{i,1}^1$
is the subgraph induced by the edge $u_1^{i,1}u_2^{i,1}$,
$P_{i,2}^1$ is the subgraph induced by the edge
$u_1^{i,2}u_2^{i,2}$, and $P_{i,3}^1$ is the subgraph induced by the
edges in $\{u_1^{i,3}u_2^{i,3},u_3^{i,3}u_4^{i,3}\}$. Furthermore,
$P_{i,2}^2$ is the subgraph induced by the edge $u_2^{i,2}u_3^{i,2}$
and $P_{i,3}^2$ is the subgraph induced by the edge
$u_2^{i,3}u_3^{i,3}$; see Figure 2 $(a)$. Note that
$F_i^1=\bigcup_{j=1}^3P_{i,j}^1$ and
$F_i^2=\bigcup_{j=2}^3P_{i,j}^2$.

Since ${\rm la}_{k}(H)=q$, it follows that there are $q$ linear
$k$-forests in $H$, say $F_1',F_2',\ldots,F_q'$. For each $F_j' \
(1\leq j\leq q)$, we assume that $Q_{j,1}',Q_{j,2}',\ldots,Q_{j,y}'$
are all the paths in $F_j'$. Then $F_j'=\bigcup_{i=1}^yQ_{j,i}'$
where $1\leq j\leq q$. Set $Q_{j,i}'=v_{j,i}^1v_{j,i}^2\ldots
v_{j,i}^b$, where $1\leq i\leq y$. For the above example, we have
$Q_{j,1}'=v_1^{j,1}v_2^{j,1}v_3^{j,1}$ and
$Q_{j,2}'=v_1^{j,2}v_2^{j,2}v_3^{j,2}v_4^{j,2}$. Then
$F_j'=Q_{j,1}'\cup Q_{j,2}'$.

We now decompose $G\times H$ into $2pq$ linear $\ell$-forests such
that each of them is formed from the paths in $F_i^1$ or $F_i^2$ and
the paths in $F_j'$, where $1\leq i\leq p$ and $1\leq j\leq q$. Note
that the subgraph $F_{i,j}^*$ induced by the edges in
$$
\{(u_1,v_1)(u_2,v_2),(u_1,v_2)(u_2,v_1)\,|\,v_1v_2\in
E(F_j'),u_1u_2\in E(F_i^1)\}
$$
and the subgraph $F_{i,j}^{**}$ induced by the edges in
$$
\{(u_1,v_1)(u_2,v_2),(u_1,v_2)(u_2,v_1)\,|\,v_1v_2\in
E(F_j'),u_1u_2\in E(F_i^2)\}
$$
are $2pq$ linear $\ell$-forests in $G\times H$, where $1\leq i\leq
p$ and $1\leq j\leq q$. Note that each edge of $G\times H$ belongs
to one of the above linear $\ell$-forests. So ${\rm
la}_{\max\{k,\ell\}}(G\times H)={\rm la}_{\ell}(G\times H)\leq 2{\rm
la}_{k}(G){\rm la}_{\ell}(H)$. \qed
\end{pf}

\subsection{For strong product}

For direct product, we have the following.
\begin{thm}\label{th2-11}
Let $G$ and $H$ be two graphs. Then
$$
\left\lceil
\frac{\Delta(G)\Delta(H)+\Delta(G)+\Delta(H)}{2}\right\rceil\\
\leq {\rm la}_{\max\{k,\ell\}}(G\boxtimes H) \leq{\rm
la}_{k}(G)+{\rm la}_{\ell}(H)+2{\rm la}_{k}(G){\rm la}_{\ell}(H).
$$
Moreover, the bounds are sharp.
\end{thm}
\begin{pf}
Note that $\Delta(G\boxtimes H)\geq
\Delta(G)\Delta(H)+\Delta(G)+\Delta(H)$. From Lemma \ref{lem2-2}, we
have
$$
{\rm la}_{\max\{k,\ell\}}(G\boxtimes H)\geq \left\lceil
\frac{\Delta(G)\Delta(H)+\Delta(G)+\Delta(H)}{2}\right\rceil.
$$
Since $E(G\boxtimes H)=E(G\times H)\cup E(G\Box H)$, it follows from
Observation \ref{obs2-2} that
\begin{eqnarray*}
{\rm la}_{\max\{k,\ell\}}(G\boxtimes H)
&\leq& {\rm la}_{\max\{k,\ell\}}(G\Box H)+{\rm la}_{\max\{k,\ell\}}(G\times H)\\
&\leq& {\rm la}_{k}(G)+{\rm la}_{\ell}(H)+2{\rm la}_{k}(G){\rm
la}_{\ell}(H),
\end{eqnarray*}
as desired.\qed
\end{pf}

\section{Results for product networks}

In this section, we demonstrate the usefulness of the proposed
constructions by applying them to some instances of Cartesian
product networks. We first study the linear $k$-arboricity of a
path, a cycle, a complete graph and a Peterson graph.
\begin{lem}\label{lem3-1}
For path $P_n \ (n\geq 2)$,
$$
\left\{
\begin{array}{ll}
{\rm la}_{k}(P_n)=1 &\mbox {\rm if} \ k\geq n-1\\
{\rm la}_{k}(P_n)=2 &\mbox {\rm if} \ 1\leq k<n-1.\\
\end{array}
\right.
$$
\end{lem}
\begin{pf}
If $k\geq n-1$, then $P_n$ itself is a line $k$-forest, and hence
${\rm la}_{k}(P_n)=1$. Suppose $k<n-1$. Let $M$ be a maximum
matching of $P_n$. Then $P_n\setminus M$ contains a matching, say
$M'$. For $1\leq k<n-1$, $M,M'$ are two line $k$-forests, and hence
${\rm la}_{k}(P_n)\leq 2$. From Observation \ref{obs2-4}, we have
${\rm la}_{k}(P_n)\geq 2$. So ${\rm la}_{k}(P_n)=2$ for $1\leq
k<n-1$.\qed
\end{pf}

\begin{lem}\label{lem3-2}
For cycle $C_n \ (n\geq 3)$,
$$
\left\{
\begin{array}{ll}
{\rm la}_{k}(C_n)=2 &\mbox {\rm if} \ n \ {\rm is \ even}, \ k\geq 1\\
{\rm la}_{k}(C_n)=3 &\mbox {\rm if} \ n \ {\rm is \ odd}, \ k=1\\
{\rm la}_{k}(C_n)=2 &\mbox {\rm if} \ n \ {\rm is \ odd}, \ k\geq 2.\\
\end{array}
\right.
$$
\end{lem}
\begin{pf}
Suppose that $n$ is even and $k\geq 1$. Since $n$ is even, it
follows that there exists a perfect matching of $G$, say $M$. Then
$M'=E(C_n)\setminus M$ is also a perfect matching of $G$. Clearly,
$M$ and $M'$ are two line $k$-forests, and hence ${\rm
la}_{k}(C_n)\leq 2$. From Observation \ref{obs2-4}, we have ${\rm
la}_{k}(C_n)=2$ for $n$ is even and $k\geq 1$.

Suppose that $n$ is odd and $k\geq 2$. Since $n$ is odd, it follows
that there exists a maximum matching of size $\frac{n-1}{2}$, say
$M$. Set $M'=E(C_n)\setminus M$. Clearly, $M$ and $M'$ are two line
$k$-forests, and hence ${\rm la}_{k}(C_n)\leq 2$. From Observation
\ref{obs2-4}, we have ${\rm la}_{k}(C_n)=2$ for $n$ is even and
$k\geq 1$.

Suppose that $n$ is odd and $k=1$. Set
$V(C_n)=\{v_1,v_2,\ldots,v_{n}\}$. We divide the edge set of $G$
into three categories: $M_1= \{v_{2i-1}v_{2i}\,|\,1\leq i\leq r\}$,
$M_2= \{v_{2i}v_{2i+2}\,|\,1\leq i\leq r\}$ and $M_3=
\{v_{1}v_{2r+1}\}$. Clearly, $M_1,M_2,M_3$ are three line
$k$-forests, and hence ${\rm la}_{k}(C_n)\leq 3$. One can easily
check that ${\rm la}_{k}(C_n)=3$ for $n$ is odd and $k=1$.\qed
\end{pf}

\begin{lem}{\upshape \cite{ChenH}}\label{lem3-3}
For complete graph $K_n \ (n\geq 2)$, ${\rm la}_{1}(K_n)=\lceil
n/2\rceil$.
\end{lem}

\begin{lem}\label{lem3-4}
For complete graph $K_n \ (n\geq 2)$, $\lceil n/2\rceil\leq {\rm
la}_{k}(K_n)\leq n$.
\end{lem}
\begin{pf}
From Lemma \ref{lem2-1}, we have
$$
\lceil n/2\rceil={\rm
la}_{n-1}(K_n)\leq {\rm la}_{k}(K_n)\leq {\rm la}_{1}(K_n)\leq
\Delta(K_n)+1=n,
$$ as desired.
\end{pf}

The Peterson graph $HP_3$ are shown in Figure 3 $(a)$. We now turn
our attention to study the linear $k$-arboricity of Peterson graphs.
\begin{lem}\label{lem3-5}
For Peterson graph graph $HP_3$, ${\rm la}_{1}(HP_3)=4$.
\end{lem}
\begin{pf}
Since $HP_3$ contains a cycle $C_5$ as its subgraph, it follows that
${\rm la}_{1}(HP_3)\geq {\rm la}_{1}(C_5)=3$. The forest $F_1$
induced by the edges in $\{v_1v_2,v_7v_{10},v_6v_9,v_3v_4\}$, the
forest $F_2$ induced by the edges in $\{v_1v_5,v_2v_3,v_8v_{10}\}$,
the forest $F_3$ induced by the edges in $\{v_4v_5,v_7v_9,v_6v_8\}$
and the forest $F_4$ induced by the edges in
$\{v_1v_6,v_2v_7,v_3v_8,v_4v_9,v_5v_{10}\}$ form $4$ linear
$1$-forests in $HP_3$. So $3\leq {\rm la}_{k}(HP_3)\leq 4$.

We claim that ${\rm la}_{k}(HP_3)=4$. Assume, to the contrary, that
${\rm la}_{k}(HP_3)=3$. Then $HP_3$ can be decomposed into $3$
linear $1$-forests, say $F_1,F_2,F_3$. Set
$C_1=v_6v_8v_{10}v_{7}v_{9}v_6$, $C_2=v_1v_2v_{3}v_{4}v_{5}v_1$, and
$M=\{v_1v_6,v_2v_7,v_3v_8,v_4v_9,v_5v_{10}\}$. We distinguish the
following cases to show this claim.

\textbf{Case 1.} $|M\cap E(F_1)|=5$ and $|M\cap E(F_2)|=|M\cap
E(F_3)|=0$.

Observe that all the edges in $C_1$ does not belong to $F_1$.
Otherwise, there is a path induced by the edges in $F_1$ such that
its length is at least $2$, which contradicts to the fact that
$k=1$. So all the edges in $C_1$ must belong to $F_2$ or $F_3$.
Since $C_1$ is a cycle of order $5$, there is a path of length at
least $2$ in $F_2$ or $F_3$, also a contradiction.

\textbf{Case 2.} $|M\cap E(F_1)|=4$, $|M\cap E(F_2)|=1$ and $|M\cap
E(F_3)|=0$.

Without loss of generality, let $v_1v_6\in F_2$ and $M\setminus
\{v_1v_6\}\subseteq E(F_1)$. Note that all the edges in $C_1$ does
not belong to $F_1$. So all the edges in $C_1$ must belong to $F_2$
or $F_3$. Since $k=1$, it follows that the elements in $E(C_1)$ must
belongs to at least $3$ linear $1$-forests, a contradiction.

\textbf{Case 3.} $|M\cap E(F_1)|=3$ and $|M\cap E(F_2)|=2$ and
$|M\cap E(F_3)|=0$, or $|M\cap E(F_1)|=3$ and $|M\cap E(F_2)|=|M\cap
E(F_3)|=1$.

From the symmetry of $HP_3$, we only need to consider the two cases
$v_1v_6,v_2v_7\notin E(F_1)$ and $v_1v_6,v_3v_8\notin E(F_1)$. At
first, we consider the former case and suppose $v_1v_6,v_2v_7\notin
E(F_1)$. Since $k=1$, it follows that $E(F_1)\cap E(C_1)=\emptyset$.
So all the edges in $C_1$ must belong to $F_2$ or $F_3$. Since
$k=1$, it follows that the elements in $E(C_1)$ must belongs to at
least $3$ linear $1$-forests, a contradiction. Next, we consider the
latter case and suppose $v_1v_6,v_3v_8\notin E(F_1)$. Since $k=1$,
it follows that $E(F_1)\cap E(C_2)=\emptyset$. So all the edges in
$C_2$ must belong to $F_2$ or $F_3$. Since $k=1$, it follows that
the elements in $E(C_2)$ must belongs to at least $3$ linear
$1$-forests, a contradiction.

\textbf{Case 4.} $|M\cap E(F_1)|=2$, $|M\cap E(F_2)|=2$ and $|M\cap
E(F_3)|=1$.

From the symmetry of $HP_3$, we have the following cases to
consider:
\begin{itemize}
\item[] (1) $v_1v_6,v_2v_7\in E(F_1)$, $v_3v_8,v_4v_9\in E(F_2)$,
$v_5v_{10}\in E(F_3)$;

\item[] (2) $v_1v_6,v_2v_7\in E(F_1)$, $v_3v_8,v_5v_{10}\in E(F_2)$,
$v_4v_{9}\in E(F_3)$;

\item[] (3) $v_1v_6,v_3v_8\in E(F_1)$, $v_2v_7,v_4v_9\in E(F_2)$, $v_5v_{10}\in
E(F_3)$;

\item[] (4) $v_1v_6,v_3v_8\in E(F_1)$, $v_2v_7,v_5v_{10}\in E(F_2)$, $v_4v_9\in
E(F_3)$.
\end{itemize}
For $(1)$, since $v_1v_6,v_2v_7\in E(F_1)$, it follows that all the
edges adjacent to $v_1v_6$ and $v_2v_7$ does not belong to $F_1$.
Then the elements in
$\{v_3v_8,v_4v_9,v_5v_{10},v_3v_{4},v_4v_{5},v_8v_{10}\}$ can belong
to $F_1$. Note that $v_3v_8,v_4v_9\in E(F_2)$, $v_5v_{10}\in
E(F_3)$, and $v_3v_{4},v_4v_{5}$ are adjacent. So $2\leq
|E(F_1)|\leq 4$. We have the following claim.

\noindent \textbf{Claim 1.} $|E(F_1)|=4$.

\noindent \textbf{Proof of Claim 1.} Assume, to the contrary, that
$|E(F_1)|=2$ or $|E(F_1)|=3$. Then $|E(F_1)\cap E(C_1)|=0$ or
$|E(F_1)\cap E(C_2)|=0$. Without loss of generality, let
$|E(F_1)\cap E(C_1)|=0$. So all the edges in $C_1$ must belong to
$F_2$ or $F_3$. Since $k=1$, it follows that the elements in
$E(C_1)$ must belongs to at least $3$ linear $1$-forests, a
contradiction.\qed\vskip 0.5mm

From Claim $1$, $|E(F_1)|=4$. Then
$F_1=\{v_1v_6,v_2v_7,v_4v_5,v_8v_{10}\}$ or $F_1=\{v_1v_6,v_2v_7,$
$v_3v_4,v_8v_{10}\}$. Suppose
$F_1=\{v_1v_6,v_2v_7,v_4v_5,v_8v_{10}\}$. Note that
$v_3v_8,v_4v_9\in E(F_2)$, $v_5v_{10}\in E(F_3)$. Then the edges in
$E(C_2)\setminus \{ v_4v_5\}$ belong to $F_2$ or $F_3$. Since the
subgraph induced by these edges is a path of length $4$, it follows
that $v_5v_1,v_2v_3\in E(F_2)$ or $v_1v_2,v_3v_4\in E(F_2)$.
Whenever which case happens, we have a path of length at least $2$
in $F_2$, a contradiction.

Similarly to the proof of $(2)$, we can also prove the correctness
of $(2)$-$(4)$.\qed
\end{pf}\vskip 0.5mm

The following observation is immediate, which will be used in lemma
\ref{lem3-6}.
\begin{obs}\label{obs3-1}
Let $C_5=w_1w_2w_3w_4w_5w_1$ be a cycle. If ${\rm la}_{3}(C_5)=2$,
then $C_5$ can be decomposed into two linear $3$-forests $F_1,F_2$
such that $F_1=w_1w_2w_3w_4,F_2=w_4w_5w_1$, or $F_1=w_1w_2w_3\cup
w_4w_5,F_2=w_3w_4\cup w_5w_1$.
\end{obs}

\begin{lem}\label{lem3-6}
For Peterson graph $HP_3$, ${\rm la}_{3}(HP_3)=3$.
\end{lem}
\begin{pf}
Note that the forest $F_1$ induced by the edges in
$\{v_2v_1v_5v_4,v_9v_7v_{10}v_8\}$, the forest $F_2$ induced by the
edges in $\{v_9v_6v_8,v_4v_3v_2\}$ and the forest $F_3$ induced by
the edges in $\{v_1v_6,v_2v_7,v_3v_8,v_4v_9,v_5v_{10}\}$ form linear
$3$-forests in $HP_3$. So ${\rm la}_{k}(HP_3)\leq 3$.

It suffices to show that ${\rm la}_{3}(HP_3)\geq 3$. Since $HP_3$
contains cycles, it follows that $2\leq {\rm la}_{k}(HP_3)\leq 3$.
We claim that ${\rm la}_{3}(HP_3)=3$. Assume, to the contrary, that
${\rm la}_{3}(HP_3)=2$. Then $HP_3$ can be decomposed into $2$
linear $3$-forests, say $F_1,F_2$. Set
$C_1=v_6v_8v_{10}v_{7}v_{9}v_6$, $C_2=v_1v_2v_{3}v_{4}v_{5}v_1$, and
$M=\{v_1v_6,v_2v_7,v_3v_8,v_4v_9,v_5v_{10}\}$. Note that $C_1,C_2$
are two cycles, $E(C_i)\cap E(F_1)\neq \emptyset$ and $E(C_i)\cap
E(F_2)\neq \emptyset$ for $i=1,2$. Since $k=3$ and ${\rm
la}_{3}(HP_3)=2$, it follows from Observation \ref{obs3-1} that
\begin{itemize}
\item[] (1) $E(C_1)\cap E(F_1)=\{v_6v_8,v_8v_{10},v_{10}v_{7}\}$ and
$E(C_1)\cap E(F_2)=\{v_7v_9,v_6v_{9}\}$.

\item[] (2) $E(C_1)\cap E(F_1)=\{v_7v_9,v_6v_{9}\}$ and $E(C_1)\cap
E(F_2)=\{v_6v_8,v_8v_{10},v_{10}v_{7}\}$.

\item[] (3) $E(C_1)\cap E(F_1)=\{v_6v_8,v_6v_{9},v_7v_{10}\}$ and $E(C_1)\cap E(F_2)=\{v_7v_9,v_8v_{10}\}$.

\item[] (4) $E(C_1)\cap E(F_1)=\{v_7v_9,v_8v_{10}\}$ and $E(C_1)\cap E(F_1)=\{v_6v_8,v_6v_{9},v_7v_{10}\}$.
\end{itemize}
By symmetry, we only need to consider $(1)$ and $(3)$. For $(1)$, we
claim that there is at most one edge in $M$ belonging to $F_1$.
Otherwise, there exists a vertex of degree $3$ in $F_1$ or a path of
length at least $4$, a contradiction. So there is at most one edge
in $M$ belonging to $F_1$. Furthermore, there are at least four
edges in $M$ belonging to $F_2$. Then there exists a path of length
at least $4$ in $F_2$, also a contradiction.

We conclude that $(3)$ holds. Similarly, $C_2$ must have the same
decomposition as $C_1$. By symmetry of $HP_3$, we have the
following.
\begin{itemize}
\item[] (3.1) $E(C_1)\cap E(F_1)=\{v_6v_8,v_6v_{9},v_7v_{10}\}$,
$E(C_1)\cap E(F_2)=\{v_7v_9,v_8v_{10}\}$, $E(C_2)\cap
E(F_1)=\{v_1v_5,v_2v_{3},v_3v_{4}\}$, $E(C_2)\cap
E(F_2)=\{v_4v_5,v_1v_{2}\}$.

\item[] (3.2) $E(C_1)\cap E(F_1)=\{v_6v_8,v_6v_{9},v_7v_{10}\}$,
$E(C_1)\cap E(F_2)=\{v_7v_9,v_8v_{10}\}$, $E(C_2)\cap
E(F_1)=\{v_1v_2,v_4v_{5}\}$, $E(C_2)\cap
E(F_2)=\{v_1v_5,v_2v_{3},v_3v_{4}\}$.

\item[] (3.3) $E(C_1)\cap E(F_1)=\{v_6v_8,v_6v_{9},v_7v_{10}\}$,
$E(C_1)\cap E(F_2)=\{v_7v_9,v_8v_{10}\}$, $E(C_2)\cap
E(F_1)=\{v_1v_2,v_1v_{5},v_3v_{4}\}$, $E(C_2)\cap
E(F_2)=\{v_4v_5,v_2v_{3}\}$.

\item[] (3.4) $E(C_1)\cap E(F_1)=\{v_6v_8,v_6v_{9},v_7v_{10}\}$,
$E(C_1)\cap E(F_2)=\{v_7v_9,v_8v_{10}\}$, $E(C_2)\cap
E(F_1)=\{v_2v_3,v_4v_{5}\}$, $E(C_2)\cap
E(F_2)=\{v_1v_2,v_1v_{5},v_3v_{4}\}$.

\item[] (3.5) $E(C_1)\cap E(F_1)=\{v_6v_8,v_6v_{9},v_7v_{10}\}$,
$E(C_1)\cap E(F_2)=\{v_7v_9,v_8v_{10}\}$, $E(C_2)\cap
E(F_1)=\{v_1v_2,v_2v_{3},v_4v_{5}\}$, $E(C_2)\cap
E(F_2)=\{v_1v_5,v_3v_{4}\}$.

\item[] (3.6) $E(C_1)\cap E(F_1)=\{v_6v_8,v_6v_{9},v_7v_{10}\}$,
$E(C_1)\cap E(F_2)=\{v_7v_9,v_8v_{10}\}$, $E(C_2)\cap
E(F_1)=\{v_1v_5,v_3v_{4}\}$, $E(C_2)\cap
E(F_2)=\{v_1v_2,v_2v_{3},v_4v_{5}\}$.
\end{itemize}
We only prove that $(3.1)$ is not true, and the other five cases can
be discussed similarly. For $(3.1)$, we claim that $v_2v_7\in
E(F_2)$. Otherwise, the path induced by the edges in
$\{v_4v_3,v_3v_{2},v_2v_{7},v_7v_{10}\}$ has length $4$, a
contradiction. So $v_2v_7\in E(F_2)$. We now focus our attention to
the edge $v_4v_9$. If $v_4v_9\in E(F_1)$, then the path induced by
the edges in $\{v_2v_3,v_3v_{4},v_4v_{9},v_9v_{6},v_6v_{8}\}$ has
length $5$, a contradiction. If $v_4v_9\in E(F_2)$, then the path
induced by the edges in
$\{v_1v_2,v_2v_{7},v_7v_{9},v_9v_{4},v_4v_{5}\}$ has length $5$, a
contradiction.\qed

From Claim 1, we have ${\rm la}_{k}(HP_3)=3$, as desired.
\end{pf}

\begin{pro}\label{pro3-1}
For a Peterson graph $HP_3$,
$$
{\rm la}_{k}(HP_3)=\left\{
\begin{array}{ll}
4 &\mbox {\rm if} \ k=1,\\
3 &\mbox {\rm if}\ k=2,\\
3 &\mbox {\rm if} \ k=3,\\
2 &\mbox {\rm if} \ k\geq 4.\\
\end{array}
\right.
$$
\end{pro}
\begin{pf}
From Lemmas \ref{lem3-4} and \ref{lem3-5}, the results follow for
$k=1,3$. For $k=2$, the forest $F_1$ induced by the paths in
$\{v_1v_5v_{10},v_6v_9v_7,v_2v_3v_8\}$, the forest $F_2$ induced by
the paths in $\{v_3v_4v_9,v_6v_8v_{10},v_1v_2v_7\}$, and the forest
$F_3$ induced by the edges in $\{v_1v_6,v_7v_{10},v_4v_5\}$  form
linear $2$-forests in $HP_3$. So ${\rm la}_{2}(HP_3)\leq 3$. From
Lemma \ref{lem2-1} and Lemma \ref{lem3-5}, we have ${\rm
la}_{2}(HP_3)\geq {\rm la}_{3}(HP_3)=3$ and hence ${\rm
la}_{2}(HP_3)=3$, as desired.
\begin{figure}[!hbpt]
\begin{center}
\includegraphics[scale=0.7]{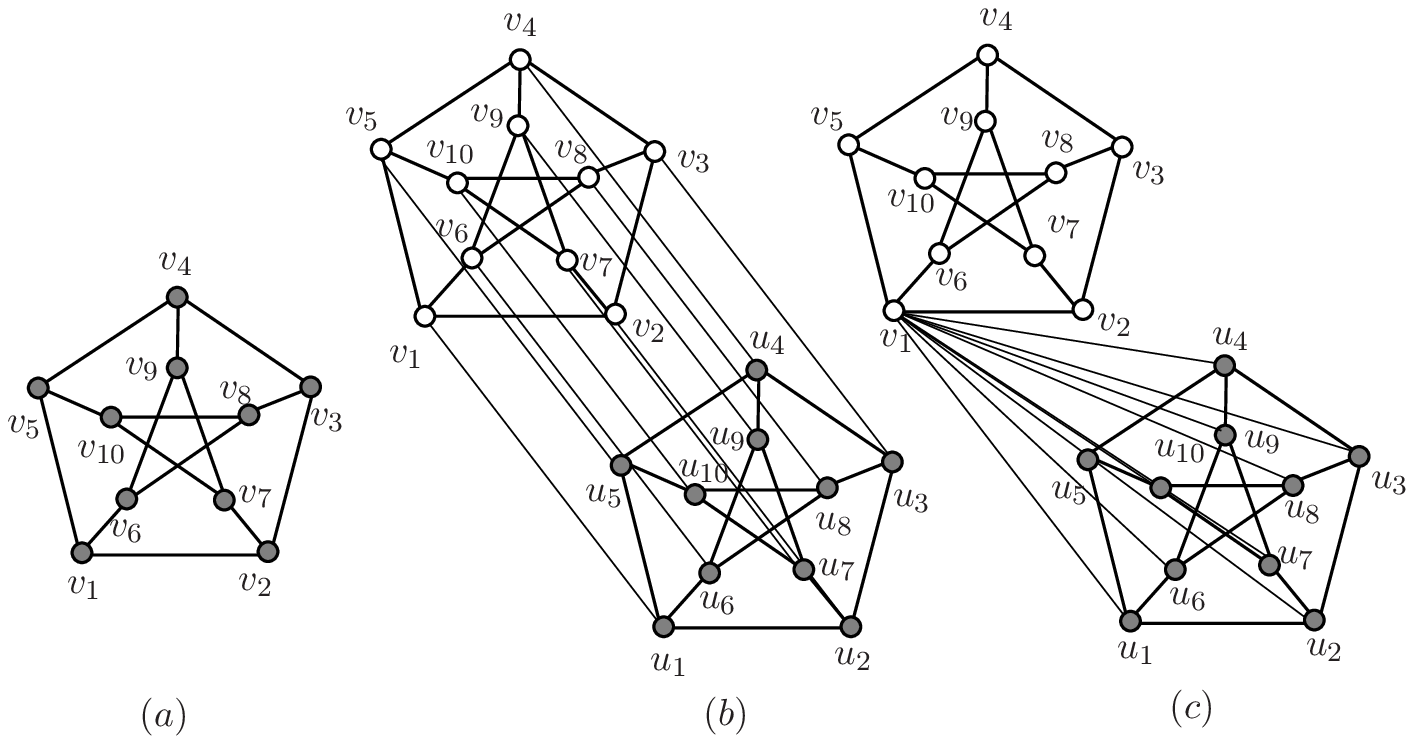}\\
Figure 3: $(a)$ Petersen graph; $(b)$ The network $HP_4$; $(c)$ The
structure of $HL_4$.
\end{center}\label{fig7}
\end{figure}
For $k\geq 4$, the forest $F_1$ induced by the paths in
$\{v_3v_2v_1v_5v_4,v_7v_9v_6v_8v_{10}\}$ and the forest $F_2$
induced by the paths in $\{v_8v_3v_4v_9,v_5v_{10}v_7v_2,v_1v_6\}$
form linear $k$-forests in $HP_3$. So ${\rm la}_{k}(HP_3)\leq 2$.
From Observation \ref{obs2-4}, we have ${\rm la}_{k}(HP_3)=2$.\qed
\end{pf}

\subsection{Two-dimensional grid graph}

A \emph{two-dimensional grid graph} is an $m\times n$ graph
$G_{n,m}$ that is the graph Cartesian product $P_n\Box P_m$ of path
graphs on $m$ and $n$ vertices. For more details on grid graph, we
refer to \cite{CalkinW, Itai}. The network $P_n\circ P_m$ is the
graph lexicographical product $P_n\circ P_m$ of path graphs on $m$
and $n$ vertices. For more details on $P_n\circ P_m$, we refer to
\cite{Mao}.
\begin{pro}\label{pro3-2}
$(i)$ For network $P_n\Box P_m \ (m\geq n\geq 3)$,
$$
\left\{
\begin{array}{ll}
{\rm la}_{k}(P_n\Box P_m)=2 &\mbox {\rm if} \ k\geq \max\{m-1,n-1\},\\
2\leq {\rm la}_{k}(P_n\Box P_m)\leq 3 &\mbox {\rm if} \ k\geq m-1, \ k\leq n-1,\\
2\leq {\rm la}_{k}(P_n\Box P_m)\leq 3 &\mbox {\rm if} \ k\geq n-1, \ k\leq m-1,\\
2\leq {\rm la}_{k}(P_n\Box P_m)\leq 4 &\mbox {\rm if} \ k\leq
\max\{m-1,n-1\}.
\end{array}
\right.
$$

$(ii)$ For network $P_n\circ P_m \ (n\geq 4, m\geq 3)$,
$$
\left\{
\begin{array}{ll}
{\rm la}_{k}(P_n\circ P_m)=m+1 &\mbox {\rm if} \ k\geq \max\{m-1,n-1\},\\
m+1\leq {\rm la}_{k}(P_n\circ P_m)\leq 2m+1 &\mbox {\rm if} \ k\geq m-1, \ k\leq n-1,\\
m+1\leq {\rm la}_{k}(P_n\circ P_m)\leq m+2 &\mbox {\rm if} \ k\geq n-1, \ k\leq m-1,\\
m+1\leq {\rm la}_{k}(P_n\circ P_m)\leq 2m+2 &\mbox {\rm if} \ k\leq
\max\{m-1,n-1\}.
\end{array}
\right.
$$

$(iii)$ For network $P_n\times P_m \ (n\geq 4, m\geq 3)$,
$$
\left\{
\begin{array}{ll}
{\rm la}_{k}(P_n\times P_m)=2 &\mbox {\rm if} \ k\geq \max\{m-1,n-1\},\\
2\leq {\rm la}_{k}(P_n\times P_m)\leq 4 &\mbox {\rm if} \ k\geq m-1, \ k\leq n-1,\\
2\leq {\rm la}_{k}(P_n\times P_m)\leq 4 &\mbox {\rm if} \ k\geq n-1, \ k\leq m-1,\\
2\leq {\rm la}_{k}(P_n\times P_m)\leq 8 &\mbox {\rm if} \ k\leq
\max\{m-1,n-1\}.
\end{array}
\right.
$$

$(iv)$ For network $P_n\boxtimes P_m \ (n\geq 4, m\geq 3)$,
$$
\left\{
\begin{array}{ll}
{\rm la}_{k}(P_n\boxtimes P_m)=4&\mbox {\rm if} \ k\geq \max\{m-1,n-1\},\\
4\leq {\rm la}_{k}(P_n\boxtimes P_m)\leq 7 &\mbox {\rm if} \ k\geq m-1, \ k\leq n-1,\\
4\leq {\rm la}_{k}(P_n\boxtimes P_m)\leq 7 &\mbox {\rm if} \ k\geq n-1, \ k\leq m-1,\\
4\leq {\rm la}_{k}(P_n\boxtimes P_m)\leq 12 &\mbox {\rm if} \ k\leq
\max\{m-1,n-1\}.
\end{array}
\right.
$$
\end{pro}
\begin{pf}
$(i)$ From Observation \ref{obs2-4}, ${\rm la}_{k}(P_n\Box P_m)\geq
2$. From Theorem \ref{th2-5}, we have
$$
{\rm la}_{k}(P_n\Box P_m)\leq {\rm la}_{k}(P_n)+{\rm
la}_{k}(P_m)\leq \left\{
\begin{array}{ll}
2 &\mbox {\rm if} \ k\geq \max\{m-1,n-1\},\\
3 &\mbox {\rm if} \ k\geq m-1, \ k\leq n-1,\\
3 &\mbox {\rm if} \ k\geq n-1, \ k\leq m-1,\\
4 &\mbox {\rm if} \ k\leq \max\{m-1,n-1\}.
\end{array}
\right.
$$
The result follows.

$(ii)$ From Theorem \ref{th2-9}, we have ${\rm la}_k(P_n\circ
P_m)\geq \left\lceil
\frac{\Delta(P_m)+|V(P_m)|\Delta(P_n)}{2}\right\rceil=m+1$. From
Theorem \ref{th2-9}, we have
$$
{\rm la}_k(P_n\circ P_m)\leq {\rm la}_{k}(P_n)|V(P_m)|+{\rm
la}_{k}(P_m)=\left\{
\begin{array}{ll}
m+1&\mbox {\rm if} \ k\geq \max\{m-1,n-1\},\\
2m+1 &\mbox {\rm if} \ k\geq m-1, \ k\leq n-1,\\
m+2 &\mbox {\rm if} \ k\geq n-1, \ k\leq m-1,\\
2m+2 &\mbox {\rm if} \ k\leq \max\{m-1,n-1\}.
\end{array}
\right.
$$

$(iii)$ From Theorem \ref{th2-10}, we have ${\rm la}_k(P_n\times
P_m)\geq \left\lceil
\frac{\Delta(P_m)\Delta(P_n)}{2}\right\rceil=2$. From Theorem
\ref{th2-10} and Lemma \ref{lem3-1}, we have
$$
{\rm la}_k(P_n\times P_m)\leq 2{\rm la}_{k}(P_n){\rm
la}_{k}(P_m)=\left\{
\begin{array}{ll}
2&\mbox {\rm if} \ k\geq \max\{m-1,n-1\},\\
4&\mbox {\rm if} \ k\geq m-1, \ k\leq n-1,\\
4&\mbox {\rm if} \ k\geq n-1, \ k\leq m-1,\\
8&\mbox {\rm if} \ k\leq \max\{m-1,n-1\}.
\end{array}
\right.
$$

$(iv)$ From Observation \ref{obs1-1}, $(i)$ and $(iii)$ of this
proposition, we have ${\rm la}_k(P_n\boxtimes P_m)\leq {\rm
la}_k(P_n\Box P_m)+{\rm la}_k(P_n\times P_m)$ and hence
$$
{\rm la}_k(P_n\boxtimes P_m)\leq \left\{
\begin{array}{ll}
4&\mbox {\rm if} \ k\geq \max\{m-1,n-1\},\\
7&\mbox {\rm if} \ k\geq m-1, \ k\leq n-1,\\
7&\mbox {\rm if} \ k\geq n-1, \ k\leq m-1,\\
12&\mbox {\rm if} \ k\leq \max\{m-1,n-1\}.
\end{array}
\right.
$$
From Theorem \ref{th2-11}, we have ${\rm la}_k(P_n\boxtimes P_m)\geq
\left\lceil
\frac{\Delta(P_m)\Delta(P_n)+\Delta(P_m)+\Delta(P_n)}{2}\right\rceil=4$.\qed
\end{pf}\\

\noindent \textbf{Remark 1.} Let $G=P_n$ and $H=P_m$. For $k\geq
\max\{m-1,n-1\}$, we have ${\rm la}_k(P_n\Box P_m)=2={\rm
la}_k(P_n)+{\rm la}_k(P_m)$, which implies that the upper bound of
Theorem \ref{th2-5} is sharp for $k=\ell$ and $k\geq
\max\{m-1,n-1\}$. For $k\geq \max\{m-1,n-1\}$, from Theorem
\ref{th2-9}, we have $m+1=\left\lceil
\frac{\Delta(P_m)+|V(P_m)|\Delta(P_n)}{2}\right\rceil\leq {\rm
la}_{k}(P_n\circ P_m)\leq {\rm la}_{k}(P_n)|V(P_m)|+{\rm
la}_{k}(P_m)=m+1$, which implies that the upper and lower bounds of
Theorem \ref{th2-9} are sharp for $k=\ell$ and $k\geq
\max\{m-1,n-1\}$. For $k\geq \max\{m-1,n-1\}$, from Theorem
\ref{th2-10}, we have $2=\lceil
\frac{\Delta(P_n)\Delta(P_m)}{2}\rceil\leq {\rm la}_{k}(P_n\times
P_m)\leq 2{\rm la}_{k}(P_n){\rm la}_{k}(P_m)=2$, which implies that
the upper and lower bounds of Theorem \ref{th2-10} are sharp for
$k=\ell$ and $k\geq \max\{m-1,n-1\}$. For $k\geq \max\{m-1,n-1\}$,
from Theorem \ref{th2-11}, we have $4=\left\lceil
\frac{\Delta(P_n)\Delta(P_m)+\Delta(P_n)+\Delta(P_m)}{2}\right\rceil
\leq {\rm la}_{k}(P_n\boxtimes P_m) \leq{\rm la}_{k}(P_n)+{\rm
la}_{k}(P_m)+2{\rm la}_{k}(P_n){\rm la}_{k}(P_m)=4$, which implies
that the upper and lower bounds of Theorem \ref{th2-11} are sharp
for $k=\ell$ and $k\geq \max\{m-1,n-1\}$.

\subsection{$r$-dimensional mesh}

An \emph{$r$-dimensional mesh} is the Cartesian product of $r$
paths. By this definition, two-dimensional grid graph is a
$2$-dimensional mesh. An $r$-dimensional hypercube is a special case
of an $r$-dimensional mesh, in which the $r$ paths are all of size
$2$; see \cite{Johnsson}.
\begin{pro}\label{pro3-2}
$(i)$ For $r$-dimensional mesh $P_{m_1}\Box P_{m_2}\Box \cdots \Box
P_{m_r}$,
$$
2\leq {\rm la}_k(P_{m_1}\Box P_{m_2}\Box \cdots \Box P_{m_r})\leq
2r.
$$

$(ii)$ For network $P_{m_1}\circ P_{m_2}\circ \cdots \circ P_{m_r}$,
$$
1+\prod_{i=2}^rm_im_{i+1}\ldots m_r\leq {\rm la}_k(P_{m_1}\circ
P_{m_2}\circ \cdots \circ P_{m_r})\leq
2\left(\sum_{i=2}^rm_im_{i+1}\ldots m_r+1\right).
$$

$(iii)$ For network $P_{m_1}\times P_{m_2}\times \cdots \times
P_{m_r}$,
$$
2^{r-1}\leq {\rm la}_k(P_{m_1}\times P_{m_2}\times \cdots \times
P_{m_r})\leq 2^{2r-1}.
$$

$(iv)$ For network $P_{m_1}\boxtimes P_{m_2}\boxtimes \cdots
\boxtimes P_{m_r}$,
$$
\frac{1}{2}(3^{r}-1)\leq {\rm la}_k(P_{m_1}\boxtimes
P_{m_2}\boxtimes \cdots \boxtimes P_{m_r})\leq 2r+2^{2r-1}.
$$
\end{pro}
\begin{pf}
$(i)$ From Corollary \ref{cor2-6} and Lemma \ref{lem3-1}, we have
${\rm la}_k(P_{m_1}\Box P_{m_2}\Box \cdots \Box P_{m_r})\leq{\rm
la}_{k}(P_{m_1})+{\rm la}_{k}(P_{m_2})+\ldots+{\rm la}_{k}(P_{m_r})
\leq 2r$. From this together with Observation \ref{obs2-4}, we have
$2\leq {\rm la}_{k}(P_{m_1})+{\rm la}_{k}(P_{m_2})+\ldots+{\rm
la}_{k}(P_{m_r}) \leq 2r$, as desired.

$(ii)$ From Theorem \ref{th2-9}, we have ${\rm la}_k(G\circ H)\leq
{\rm la}_{k}(G)|V(H)|+{\rm la}_{k}(H)$ for any two graphs $G$ and
$H$, and hence
\begin{eqnarray*} \label{n-t=0}
&&{\rm la}_k(P_{m_1}\circ P_{m_2}\circ \cdots \circ P_{m_r})\\
&=&{\rm la}_k((P_{m_1}\circ P_{m_2}\circ \cdots \circ P_{m_{r-1}})\circ P_{m_r})\\
&\leq&{\rm la}_k(P_{m_1}\circ P_{m_2}\circ \cdots \circ P_{m_{r-1}})m_r+{\rm la}_k(P_{m_r})\\
&\leq&{\rm la}_k(P_{m_1}\circ P_{m_2}\circ \cdots \circ P_{m_{r-1}})m_r+2\\
&\leq&[{\rm la}_k(P_{m_1}\circ P_{m_2}\circ \cdots \circ P_{m_{r-2}})m_{r-1}+{\rm la}_k(P_{m_{r-1}})]m_r+2\\
&\leq&[{\rm la}_k(P_{m_1}\circ P_{m_2}\circ \cdots \circ P_{m_{r-2}})m_{r-1}+2]m_r+2\\
&=&{\rm la}_k(P_{m_1}\circ P_{m_2}\circ \cdots \circ P_{m_{r-2}})m_{r-1}m_r+2m_r+2\\
&\leq&[{\rm la}_k(P_{m_1}\circ P_{m_2}\circ \cdots \circ P_{m_{r-3}})m_{r-2}+{\rm la}_k(P_{m_{r-2}})]m_{r-1}m_r+2m_r+2\\
&\leq&{\rm la}_k(P_{m_1}\circ P_{m_2}\circ \cdots \circ P_{m_{r-3}})m_{r-2}m_{r-1}m_r+2m_{r-1}m_r+2m_r+2\\
&\leq&\ldots\\
&\leq&{\rm la}_k(P_{m_1}\circ P_{m_2})m_3\ldots
m_{r-2}m_{r-1}m_r+2m_4\ldots m_{r-2}m_{r-1}m_r+\ldots+2m_{r-1}m_r+2m_r+2\\
&\leq&2m_2m_3\ldots
m_{r-2}m_{r-1}m_r+2m_3\ldots m_{r-2}m_{r-1}m_r+\ldots+2m_{r-1}m_r+2m_r+2\\
&\leq&2\sum_{i=2}^rm_im_{i+1}\ldots m_r+2.
\end{eqnarray*}
From Theorem \ref{th2-9}, we have ${\rm la}_k(G\circ H)\geq
\left\lceil \frac{\Delta(H)+|V(H)|\Delta(G)}{2}\right\rceil$ for any
two graphs $G$ and $H$, and hence
\begin{eqnarray*} \label{n-t=0}
&&{\rm la}_k(P_{m_1}\circ P_{m_2}\circ \cdots \circ P_{m_r})\\
&\geq&\left\lceil\frac{\Delta(P_{m_r})+|V(P_{m_r})|\Delta(P_{m_1}\circ P_{m_2}\circ \cdots \circ P_{m_{r-1}})}{2}\right\rceil\\
&\geq&1+\left\lceil \frac{m_r}{2}\Delta(P_{m_1}\circ P_{m_2}\circ \cdots \circ P_{m_{r-1}})\right\rceil\\
&\geq&1+\left\lceil \frac{m_r}{2}\left(\Delta(P_{m_{r-1}})+|V(P_{m_{r-1}})|\Delta(P_{m_1}\circ P_{m_2}\circ \cdots \circ P_{m_{r-2}})\right)\right\rceil\\
&\geq&1+m_r+\left\lceil \frac{m_rm_{r-1}}{2}\Delta(P_{m_1}\circ P_{m_2}\circ \cdots \circ P_{m_{r-2}})\right\rceil\\
&\geq&\ldots\\
&\geq&1+m_r+m_rm_{r-1}+\cdots +\left\lceil \frac{m_rm_{r-1}\cdots m_2}{2}\Delta(P_{m_1})\right\rceil\\
&\geq&1+m_r+m_rm_{r-1}+\cdots +m_rm_{r-1}\cdots m_2\\
&\geq&1+\prod_{i=2}^rm_im_{i+1}\ldots m_r.
\end{eqnarray*}

$(iii)$ From Theorem \ref{th2-10}, we have ${\rm la}_{k}(G\times
H)\leq 2{\rm la}_{k}(G){\rm la}_{k}(H)$ for any two graphs $G$ and
$H$, and hence
\begin{eqnarray*} \label{n-t=0}
{\rm la}_k(P_{m_1}\times P_{m_2}\times \cdots \times P_{m_r})
&\leq&2{\rm la}_k(P_{m_1}\times P_{m_2}\times \cdots \times P_{m_{r-1}}){\rm la}_k(P_{m_r})\\
&\leq&2^2{\rm la}_k(P_{m_1}\times P_{m_2}\times \cdots \times P_{m_{r-1}})\\
&\leq&2^3{\rm la}_k(P_{m_1}\times P_{m_2}\times \cdots \times P_{m_{r-2}}){\rm la}_k(P_{m_{r-1}})\\
&\leq&2^4{\rm la}_k(P_{m_1}\times P_{m_2}\times \cdots \times P_{m_{r-2}})\\
&\leq&\ldots\\
&\leq&2^{2(r-1)}{\rm la}_k(P_{m_1})\\
&\leq&2^{2r-1}.
\end{eqnarray*}
From Theorem \ref{th2-10}, we have ${\rm la}_k(G\times H)\geq
\left\lceil \frac{\Delta(G)\Delta(H)}{2}\right\rceil$ for any two
graphs $G$ and $H$, and hence
\begin{eqnarray*} \label{n-t=0}
{\rm la}_k(P_{m_1}\times P_{m_2}\times \cdots \times P_{m_r})
&\geq&\left\lceil\frac{\Delta(P_{m_1}\times P_{m_2}\times \cdots \times P_{m_{r-1}})\Delta(P_{m_r})}{2}\right\rceil\\
&\geq&\Delta(P_{m_1}\times P_{m_2}\times \cdots \times P_{m_{r-1}})\\
&\geq&\Delta(P_{m_1}\times P_{m_2}\times \cdots \times P_{m_{r-2}})\Delta(P_{m_{r-1}})\\
&=&2\Delta(P_{m_1}\times P_{m_2}\times \cdots \times P_{m_{r-2}})\\
&\geq&2\Delta(P_{m_1}\times P_{m_2}\times \cdots \times P_{m_{r-3}})\Delta(P_{m_{r-2}})\\
&=&2^2\Delta(P_{m_1}\times P_{m_2}\times \cdots \times P_{m_{r-3}})\\
&\geq&\ldots\\
&\geq&2^{r-2}{\rm la}_k(P_{m_1})\\
&\geq&2^{r-1}.
\end{eqnarray*}

$(iv)$ From $(i),(iii)$ of this proposition and Observation
\ref{obs1-1}, we have
\begin{eqnarray*} \label{n-t=0}
{\rm la}_k(P_{m_1}\boxtimes P_{m_2}\boxtimes \cdots \boxtimes
P_{m_r}) &\leq&{\rm la}_k(P_{m_1}\Box P_{m_2}\Box \cdots \Box
P_{m_r})+{\rm la}_k(P_{m_1}\times P_{m_2}\times \cdots \times
P_{m_r})\\
&\leq&2r+2^{2r-1}.
\end{eqnarray*}
From Theorem \ref{th2-10}, we have
$$
{\rm la}_{k}(G\boxtimes H)\geq \left\lceil
\frac{\Delta(G)\Delta(H)+\Delta(G)+\Delta(H)}{2}\right\rceil=\left\lceil
\frac{\Delta(G)(\Delta(H)+1)+\Delta(H)}{2}\right\rceil
$$
for
any two graphs $G$ and $H$, and hence
\begin{eqnarray*} \label{n-t=0}
&&{\rm la}_k(P_{m_1}\boxtimes P_{m_2}\boxtimes \cdots \boxtimes
P_{m_r})\\
&\geq&\left\lceil \frac{\Delta(P_{m_1}\boxtimes P_{m_2}\boxtimes
\cdots \boxtimes
P_{m_{r-1}})(\Delta(P_{m_r})+1)+\Delta(P_{m_r})}{2}\right\rceil\\
&=&\left\lceil \frac{3}{2}\Delta(P_{m_1}\boxtimes P_{m_2}\boxtimes
\cdots \boxtimes
P_{m_{r-1}})\right\rceil+1\\
&\geq&\left\lceil \frac{3}{2}\Delta(P_{m_1}\boxtimes
P_{m_2}\boxtimes \cdots \boxtimes
P_{m_{r-2}})(\Delta(P_{m_{r-1}})+1)+\Delta(P_{m_{r-1}})\right\rceil+1\\
&=&\left\lceil \frac{3^2}{2}\Delta(P_{m_1}\boxtimes P_{m_2}\boxtimes
\cdots \boxtimes
P_{m_{r-2}})\right\rceil+3+1\\
&\geq&\left\lceil \frac{3^2}{2}\Delta(P_{m_1}\boxtimes
P_{m_2}\boxtimes \cdots \boxtimes
P_{m_{r-3}})(\Delta(P_{m_{r-2}})+1)+\Delta(P_{m_{r-2}})\right\rceil+1\\
&=&\left\lceil \frac{3^3}{2}\Delta(P_{m_1}\boxtimes P_{m_2}\boxtimes
\cdots \boxtimes
P_{m_{r-3}})\right\rceil+3^2+3+1\\
&\geq&\ldots\\
&\geq&\left\lceil \frac{3^{r-1}}{2}\Delta(P_{m_1})\right\rceil+3^{r-2}+\ldots+3^2+3+1\\
&=&\frac{1}{2}(3^{r}-1).
\end{eqnarray*}
\end{pf}

\subsection{$r$-dimensional torus}

An \emph{$r$-dimensional torus} is the Cartesian product of $r$
cycles $C_{m_1},C_{m_2},\cdots,C_{m_r}$ of size at least three. The
cycles $C_{m_i}$ are not necessary to have the same size. Ku et al.
\cite{Ku} showed that there are $r$ edge-disjoint spanning trees in
an $r$-dimensional torus. The network $C_{m_1}\circ C_{m_2}\circ
\cdots \circ C_{m_r}$ is investigated in \cite{Mao}. Here, we
consider the networks constructed by $C_{m_1}\Box C_{m_2}\Box \cdots
\Box C_{m_r}$ and $C_{m_1}\circ C_{m_2}\circ \cdots \circ C_{m_r}$,
respectively.
\begin{pro}\label{pro3-3}
$(i)$ For $r$-dimensional torus $C_{m_1}\Box C_{m_2}\Box \cdots \Box
C_{m_r}$,
$$
2\leq {\rm la}_k(C_{m_1}\Box C_{m_2}\Box \cdots \Box C_{m_r})\leq
3r.
$$

$(ii)$ For network $C_{m_1}\circ C_{m_2}\circ \cdots \circ C_{m_r}$,
$$
1+\prod_{i=2}^rm_im_{i+1}\ldots m_r\leq {\rm la}_k(C_{m_1}\circ
C_{m_2}\circ \cdots \circ C_{m_r})\leq
3\left(\sum_{i=2}^rm_im_{i+1}\ldots m_r+1\right).
$$

$(ii)$ For network $C_{m_1}\times C_{m_2}\times \cdots \times
C_{m_r}$,
$$
2^{r-1}\leq {\rm la}_k(C_{m_1}\times C_{m_2}\times \cdots \times
C_{m_r})\leq 3\cdot 6^{r-1}.
$$

$(iv)$ For network $C_{m_1}\boxtimes C_{m_2}\boxtimes \cdots
\boxtimes C_{m_r}$,
$$
\frac{1}{2}(3^{r}-1)\leq {\rm la}_k(C_{m_1}\boxtimes
C_{m_2}\boxtimes \cdots \boxtimes C_{m_r})\leq 3(r+6^{r-1}).
$$
\end{pro}
\begin{pf}
$(i)$ From Corollary \ref{cor2-6} and Lemma \ref{lem3-4}, we have
\begin{eqnarray*}
2&\leq &\max\{{\rm la}_{k}(C_{m_1}),{\rm
la}_{k}(C_{m_2}),\ldots,{\rm
la}_{k}(C_{m_r})\}\\
&\leq& {\rm la}_k(C_{m_1}\Box C_{m_2}\Box \cdots \Box C_{m_r})\\
&\leq&{\rm la}_{k}(C_{m_1})+{\rm la}_{k}(C_{m_2})+\ldots+{\rm
la}_{k}(C_{m_r})\\
&\leq&3r.
\end{eqnarray*}

$(ii)$ From Theorem \ref{th2-5}, we have ${\rm la}_k(G\circ H)\leq
{\rm la}_{k}(G)|V(H)|+{\rm la}_{k}(H)$ for any two graphs $G$ and
$H$, and hence
\begin{eqnarray*} \label{n-t=0}
&&{\rm la}_k(C_{m_1}\circ C_{m_2}\circ \cdots \circ C_{m_r})\\
&=&{\rm la}_k((C_{m_1}\circ C_{m_2}\circ \cdots \circ C_{m_{r-1}})\circ C_{m_r})\\
&\leq&{\rm la}_k(C_{m_1}\circ C_{m_2}\circ \cdots \circ C_{m_{r-1}})m_r+{\rm la}_k(C_{m_r})\\
&\leq&{\rm la}_k(C_{m_1}\circ C_{m_2}\circ \cdots \circ C_{m_{r-1}})m_r+3\\
&\leq&[{\rm la}_k(C_{m_1}\circ C_{m_2}\circ \cdots \circ C_{m_{r-2}})m_{r-1}+{\rm la}_k(C_{m_{r-1}})]m_r+3\\
&\leq&[{\rm la}_k(C_{m_1}\circ C_{m_2}\circ \cdots \circ C_{m_{r-2}})m_{r-1}+3]m_r+3\\
&=&{\rm la}_k(C_{m_1}\circ C_{m_2}\circ \cdots \circ C_{m_{r-2}})m_{r-1}m_r+3m_r+3\\
&\leq&[{\rm la}_k(C_{m_1}\circ C_{m_2}\circ \cdots \circ C_{m_{r-3}})m_{r-2}+{\rm la}_k(C_{m_{r-2}})]m_{r-1}m_r+3m_r+3\\
&\leq&{\rm la}_k(C_{m_1}\circ C_{m_2}\circ \cdots \circ C_{m_{r-3}})m_{r-2}m_{r-1}m_r+3m_{r-1}m_r+3m_r+3\\
&\leq&\ldots\\
&\leq&{\rm la}_k(C_{m_1}\circ C_{m_2})m_3\ldots
m_{r-2}m_{r-1}m_r+3m_4\ldots m_{r-2}m_{r-1}m_r+\ldots+3m_{r-1}m_r+3m_r+3\\
&\leq&3m_2m_3\ldots
m_{r-2}m_{r-1}m_r+3m_3\ldots m_{r-2}m_{r-1}m_r+\ldots+3m_{r-1}m_r+3m_r+3\\
&\leq&3\sum_{i=2}^rm_im_{i+1}\ldots m_r+3.
\end{eqnarray*}
From Observation \ref{obs2-1} and $(ii)$ of Proposition
\ref{pro3-2}, we have ${\rm la}_k(C_{m_1}\circ C_{m_2}\circ \cdots
\circ C_{m_r})\geq {\rm la}_k(P_{m_1}\circ P_{m_2}\circ \cdots \circ
P_{m_r})\geq 1+\prod_{i=2}^rm_im_{i+1}\ldots m_r$.

$(iii)$ From Theorem \ref{th2-10}, we have ${\rm la}_{k}(G\times
H)\leq 2{\rm la}_{k}(G){\rm la}_{k}(H)$ for any two graphs $G$ and
$H$, and hence
\begin{eqnarray*} \label{n-t=0}
{\rm la}_k(C_{m_1}\times C_{m_2}\times \cdots \times C_{m_r})
&\leq&2{\rm la}_k(C_{m_1}\times C_{m_2}\times \cdots \times C_{m_{r-1}}){\rm la}_k(C_{m_r})\\
&\leq&6{\rm la}_k(C_{m_1}\times C_{m_2}\times \cdots \times C_{m_{r-1}})\\
&\leq&6[2{\rm la}_k(C_{m_1}\times C_{m_2}\times \cdots \times C_{m_{r-2}}){\rm la}_k(C_{m_{r-1}})]\\
&\leq&6^2{\rm la}_k(C_{m_1}\times C_{m_2}\times \cdots \times C_{m_{r-2}})\\
&\leq&\ldots\\
&\leq&6^{r-1}{\rm la}_k(C_{m_1})\\
&\leq&3\cdot 6^{r-1}.
\end{eqnarray*}
From Observation \ref{obs2-1} and $(iii)$ of Proposition
\ref{pro3-2}, we have ${\rm la}_k(C_{m_1}\times C_{m_2}\times \cdots
\times C_{m_r})\geq {\rm la}_k(P_{m_1}\times P_{m_2}\times \cdots
\times P_{m_r})\geq 2^{r-1}$.

$(iv)$ From $(i),(iv)$ of this proposition and Observation
\ref{obs1-1}, we have
\begin{eqnarray*} \label{n-t=0}
&&{\rm la}_k(C_{m_1}\boxtimes C_{m_2}\boxtimes \cdots \boxtimes
C_{m_r})\\
&\leq&{\rm la}_k(C_{m_1}\Box C_{m_2}\Box \cdots \Box C_{m_r})+{\rm
la}_k(C_{m_1}\times C_{m_2}\times \cdots \times
C_{m_r})\\
&\leq&3(r+6^{r-1}).
\end{eqnarray*}
From Observation \ref{obs2-1} and $(iv)$ of Proposition
\ref{pro3-2}, we have ${\rm la}_k(C_{m_1}\boxtimes C_{m_2}\boxtimes
\cdots \boxtimes C_{m_r})\geq {\rm la}_k(P_{m_1}\boxtimes
P_{m_2}\boxtimes \cdots \boxtimes P_{m_r})\geq
1+\prod_{i=2}^rm_im_{i+1}\ldots m_r$.

\end{pf}

\subsection{$r$-dimensional generalized hypercube}

Let $K_m$ be a clique of $m$ vertices, $m\geq 2$. An
\emph{$r$-dimensional generalized hypercube} \cite{DayA,
Fragopoulou} is the Cartesian product of $r$ cliques. We have the
following:
\begin{pro}\label{pro3-5}
$(i)$ For generalized hypercube $K_{m_1}\Box K_{m_2}\Box \cdots \Box
K_{m_r} \ (m_i\geq 2, \ r\geq 2, \ 1\leq i\leq r)$,
$$
\max\left\{\left\lceil \frac{m_i}{2}\right\rceil\,|\,1\leq i\leq
r\right\}\leq {\rm la}_k(K_{m_1}\Box K_{m_2}\Box \cdots \Box
K_{m_r})\leq \sum_{i=1}^r m_i.
$$

$(ii)$ For network $K_{m_1}\circ K_{m_2}\circ \cdots \circ K_{m_r} \
(m_i\geq 2, \ r\geq 2, \ 1\leq i\leq r)$,
$$
\left\lceil\frac{\sum_{i=1}^r m_i}{2}\right\rceil\leq {\rm
la}_k(K_{m_1}\Box K_{m_2}\Box \cdots \Box K_{m_r})\leq
\frac{\sum_{i=1}^r m_i}{2}.
$$

$(iii)$ For network $K_{m_1}\times K_{m_2}\times \cdots \times
K_{m_r} \ (m_i\geq 2, \ r\geq 2, \ 1\leq i\leq r)$,
$$
\left\lceil \frac{1}{2}\prod_{i=1}^r(m_i-1)\right\rceil\leq {\rm
la}_k(K_{m_1}\times K_{m_2}\times \cdots \times K_{m_r})\leq
2^{r-1}\prod_{i=1}^rm_i.
$$

$(iv)$ For network $K_{m_1}\boxtimes K_{m_2}\boxtimes \cdots
\boxtimes K_{m_r} \ (m_i\geq 2, \ r\geq 2, \ 1\leq i\leq r)$,
$$
\left\lceil \frac{1}{2}\prod_{i=1}^rm_rm_{r-1}\ldots
m_{i+1}(m_{i}-1)\right\rceil\leq {\rm la}_k(K_{m_1}\boxtimes
K_{m_2}\boxtimes \cdots \boxtimes K_{m_r})\leq \sum_{i=1}^r
m_i+2^{r-1}\prod_{i=1}^rm_i.
$$
\end{pro}
\begin{pf}
$(i)$ From Corollary \ref{cor2-6} and Lemma \ref{lem3-4}, we have
\begin{eqnarray*}
\max\left\{\left\lceil \frac{m_i}{2}\right\rceil\,|\,1\leq i\leq
r\right\}&=&\max\{{\rm la}_{k}(K_{m_1}),{\rm
la}_{k}(K_{m_2}),\ldots,{\rm
la}_{k}(K_{m_r})\}\\
&\leq& {\rm la}_k(K_{m_1}\Box K_{m_2}\Box \cdots \Box K_{m_r})\\
&\leq&{\rm la}_{k}(K_{m_1})+{\rm la}_{k}(K_{m_2})+\ldots+{\rm
la}_{k}(K_{m_r})\\
&\leq& \sum_{i=1}^r m_i.
\end{eqnarray*}

$(ii)$ From the definition of lexicographical product, $K_{m_1}\circ
K_{m_2}\circ \cdots \circ K_{m_r}$ is a complete graph. From Lemma
\ref{lem3-4}, we have
$$
\left\lceil\frac{\sum_{i=1}^r m_i}{2}\right\rceil\leq {\rm
la}_k(K_{m_1}\circ K_{m_2}\circ \cdots \circ K_{m_r})\leq
\frac{\sum_{i=1}^r m_i}{2}.
$$

$(iii)$ From Theorem \ref{th2-10}, we have ${\rm la}_{k}(G\times
H)\leq 2{\rm la}_{k}(G){\rm la}_{k}(H)$ for any two graphs $G$ and
$H$, and hence
\begin{eqnarray*} \label{n-t=0}
{\rm la}_{k}(K_{m_1}\times K_{m_2}\times \cdots \times K_{m_r})
&\leq&2{\rm la}_{k}(K_{m_1}\times K_{m_2}\times \cdots \times K_{m_{r-1}}){\rm la}_{k}(K_{m_r})\\
&\leq&2m_r{\rm la}_{k}(K_{m_1}\times K_{m_2}\times \cdots \times K_{m_{r-1}})\\
&\leq&2^2m_r{\rm la}_{k}(K_{m_1}\times K_{m_2}\times \cdots \times K_{m_{r-2}}){\rm la}_{k}(K_{m_{r-1}})\\
&\leq&2^2m_rm_{r-1}{\rm la}_{k}(K_{m_1}\times K_{m_2}\times \cdots \times K_{m_{r-2}})\\
&\leq&\ldots\\
&\leq&2^{r-1}m_rm_{r-1}\ldots m_2{\rm la}_{k}(K_{m_1})\\
&=&2^{r-1}\prod_{i=1}^rm_i.
\end{eqnarray*}
From Theorem \ref{th2-10}, we have ${\rm la}_{k}(G\times H)\geq
\left\lceil \frac{\Delta(G)\Delta(H)}{2}\right\rceil$ for any two
graphs $G$ and $H$, and hence
\begin{eqnarray*} \label{n-t=0}
{\rm la}_{k}(K_{m_1}\times K_{m_2}\times \cdots \times K_{m_r})
&\geq&\left\lceil \frac{\Delta(K_{m_1}\times K_{m_2}\times \cdots \times K_{m_{r-1}})\Delta(K_{m_r})}{2}\right\rceil\\
&\geq&\left\lceil \frac{m_r-1}{2}\Delta(K_{m_1}\times K_{m_2}\times \cdots \times K_{m_{r-1}})\right\rceil\\
&\geq&\left\lceil \frac{m_r-1}{2}\Delta(K_{m_1}\times K_{m_2}\times \cdots \times K_{m_{r-2}})\Delta(K_{m_r})\right\rceil\\
&\geq&\left\lceil \frac{(m_r-1)(m_{r-1}-1)}{2}\Delta(K_{m_1}\times K_{m_2}\times \cdots \times K_{m_{r-2}})\right\rceil\\
&\geq&\ldots\\
&\geq&\left\lceil \frac{(m_r-1)(m_{r-1}-1)\ldots (m_{2}-1)}{2}\Delta(K_{m_1})\right\rceil\\
&=&\left\lceil \frac{1}{2}\prod_{i=1}^r(m_i-1)\right\rceil.
\end{eqnarray*}

$(iv)$ From $(i),(iii)$ of this proposition and Observation
\ref{obs1-1}, we have
\begin{eqnarray*} \label{n-t=0}
&&{\rm la}_k(K_{m_1}\boxtimes K_{m_2}\boxtimes \cdots \boxtimes
K_{m_r})\\
&\leq&{\rm la}_k(K_{m_1}\Box K_{m_2}\Box \cdots \Box K_{m_r})+{\rm
la}_k(K_{m_1}\times K_{m_2}\times \cdots \times
K_{m_r})\\
&\leq&\sum_{i=1}^r m_i+2^{r-1}\prod_{i=1}^rm_i.
\end{eqnarray*}
From Theorem \ref{th2-11}, we have
$$
{\rm la}_{k}(G\boxtimes H)\geq \left\lceil
\frac{\Delta(G)(\Delta(H)+1)+\Delta(H)}{2}\right\rceil
$$
for any two graphs $G$ and $H$, and hence
\begin{eqnarray*} \label{n-t=0}
&&{\rm la}_k(K_{m_1}\boxtimes K_{m_2}\boxtimes \cdots \boxtimes
K_{m_r})\\
&\geq&\left\lceil \frac{\Delta(K_{m_1}\boxtimes K_{m_2}\boxtimes
\cdots \boxtimes
K_{m_{r-1}})(\Delta(K_{m_r})+1)+\Delta(K_{m_r})}{2}\right\rceil\\
&=&\left\lceil \frac{m_r-1}{2}+\frac{m_r}{2}\Delta(K_{m_1}\boxtimes
K_{m_2}\boxtimes \cdots \boxtimes
K_{m_{r-1}})\right\rceil\\
&\geq&\left\lceil
\frac{m_r-1}{2}+\frac{m_r}{2}\left[\Delta(K_{m_1}\boxtimes
K_{m_2}\boxtimes \cdots \boxtimes
K_{m_{r-2}})(\Delta(K_{m_{r-1}})+1)+\Delta(K_{m_{r-1}})\right]\right\rceil\\
&=&\left\lceil
\frac{m_r-1}{2}+\frac{m_r(m_{r-1}-1)}{2}+\frac{m_rm_{r-1}}{2}\Delta(K_{m_1}\boxtimes
K_{m_2}\boxtimes \cdots \boxtimes
K_{m_{r-2}})\right\rceil\\
&\geq&\ldots\\
&\geq&\left\lceil \frac{m_r-1}{2}+\frac{m_r(m_{r-1}-1)}{2}+\frac{m_rm_{r-1}(m_{r-2}-1)}{2}+\ldots+\frac{m_rm_{r-1}\ldots m_{2}}{2}\Delta(K_{m_1})\right\rceil\\
&=&\left\lceil \frac{1}{2}\prod_{i=1}^rm_rm_{r-1}\ldots
m_{i+1}(m_{i}-1)\right\rceil.
\end{eqnarray*}\qed
\end{pf}

\subsection{$n$-dimensional hyper Petersen network}

An \emph{$n$-dimensional hyper Petersen network} $HP_n$ is the
product of the well-known Petersen graph and $Q_{n-3}$ \cite{Das},
where $n\geq 3$ and $Q_{n-3}$ denotes an $(n-3)$-dimensional
hypercube. The cases $n=3$ and $4$ of hyper Petersen networks are
depicted in Figure 2. Note that $HP_3$ is just the Petersen graph;
see Figure 3 $(a)$.

The network $HL_n$ is the lexicographical product of the Petersen
graph and $Q_{n-3}$, where $n\geq 3$ and $Q_{n-3}$ denotes an
$(n-3)$-dimensional hypercube; see \cite{Mao}. Note that $HL_4$ is a
graph obtained from two copies of the Petersen graph by add one edge
between one vertex in a copy of the Petersen graph and one vertex in
another copy; see Figure 3 $(c)$ for an example (We only show the
edges $v_1u_i \ (1\leq i\leq 10)$).

Similarly, the networks $HD_n$ and $HS_n$ are defined as the direct
and strong product of the Petersen graph and $Q_{n-3}$,
respectively, where $n\geq 3$ and $Q_{n-3}$ denotes an
$(n-3)$-dimensional hypercube. Note that $HL_3=HD_3=HS_3$ is just
the Petersen graph, and
\begin{pro}\label{pro3-5}
$(i)$ For network $HP_4$,
$$
\left\{
\begin{array}{ll}
4\leq {\rm la}_{k}(HP_4)\leq 5 &\mbox {\rm if} \ k=1,\\
3\leq {\rm la}_{k}(HP_4)\leq 4 &\mbox {\rm if}\ k=2,\\
3\leq {\rm la}_{k}(HP_4)\leq 4 &\mbox {\rm if} \ k=3,\\
2\leq {\rm la}_{k}(HP_4)\leq 3 &\mbox {\rm if} \ k\geq 4.\\
\end{array}
\right.
$$

$(ii)$ For network $HL_4$,
$$
\left\{
\begin{array}{ll}
4\leq {\rm la}_{k}(HL_4)\leq 14 &\mbox {\rm if} \ k=1,\\
4\leq {\rm la}_{k}(HL_4)\leq 13 &\mbox {\rm if}\ k=2,\\
4\leq {\rm la}_{k}(HL_4)\leq 13 &\mbox {\rm if} \ k=3,\\
4\leq {\rm la}_{k}(HL_4)=12 &\mbox {\rm if} \ k\geq 4.\\
\end{array}
\right.
$$

$(iii)$ For network $HD_4$,
$$
\left\{
\begin{array}{ll}
2\leq {\rm la}_{k}(HD_4)\leq 8 &\mbox {\rm if} \ k=1,\\
2\leq {\rm la}_{k}(HD_4)\leq 6 &\mbox {\rm if}\ k=2,\\
2\leq {\rm la}_{k}(HD_4)\leq 6 &\mbox {\rm if} \ k=3,\\
2\leq {\rm la}_{k}(HD_4)\leq 4 &\mbox {\rm if} \ k\geq 4.\\
\end{array}
\right.
$$

$(iv)$ For network $HS_4$,
$$
\left\{
\begin{array}{ll}
4\leq {\rm la}_{k}(HS_4)\leq 13 &\mbox {\rm if} \ k=1,\\
4\leq {\rm la}_{k}(HS_4)\leq 10 &\mbox {\rm if}\ k=2,\\
4\leq {\rm la}_{k}(HS_4)\leq 10 &\mbox {\rm if} \ k=3,\\
4\leq {\rm la}_{k}(HS_4)\leq 7 &\mbox {\rm if} \ k\geq 4.\\
\end{array}
\right.
$$
\end{pro}
\begin{pf}
$(i)$ From Theorem \ref{th2-5}, we have
$$
{\rm la}_{k}(HP_4)\geq {\rm la}_{k}(HP_3)\geq \left\{
\begin{array}{ll}
4 &\mbox {\rm if} \ k=1,\\
3 &\mbox {\rm if}\ k=2,\\
3 &\mbox {\rm if} \ k=3,\\
2 &\mbox {\rm if} \ k\geq 4.\\
\end{array}
\right.
$$
From Theorem \ref{th2-5}, Lemmas \ref{lem3-1} and \ref{lem3-5}, we
have
$$
{\rm la}_{k}(HP_4)\leq {\rm la}_{k}(HP_3)+{\rm la}_{k}(P_2)={\rm
la}_{k}(HP_3)+1\leq \left\{
\begin{array}{ll}
5 &\mbox {\rm if} \ k=1,\\
4 &\mbox {\rm if}\ k=2,\\
4 &\mbox {\rm if} \ k=3,\\
3 &\mbox {\rm if} \ k\geq 4.\\
\end{array}
\right.
$$

$(ii)$ From Theorem \ref{th2-9}, Lemmas \ref{lem3-1} and
\ref{lem3-5}, we have
$$
4\leq {\rm la}_{k}(HP_4)\leq {\rm la}_{k}(HP_3)|V(P_2)|+{\rm
la}_{k}(P_2)=2{\rm la}_{k}(HP_3)+1\leq \left\{
\begin{array}{ll}
9 &\mbox {\rm if} \ k=1,\\
7 &\mbox {\rm if}\ k=2,\\
7 &\mbox {\rm if} \ k=3,\\
5 &\mbox {\rm if} \ k\geq 4.\\
\end{array}
\right.
$$

$(iii)$ From Theorem \ref{th2-10}, Lemmas \ref{lem3-1} and
\ref{lem3-5}, we have
$$
{\rm la}_{k}(HD_4)\leq 2{\rm la}_{k}(HP_3){\rm la}_{k}(P_2)=2{\rm
la}_{k}(HP_3)\leq \left\{
\begin{array}{ll}
8 &\mbox {\rm if} \ k=1,\\
6 &\mbox {\rm if}\ k=2,\\
6 &\mbox {\rm if} \ k=3,\\
4 &\mbox {\rm if} \ k\geq 4.\\
\end{array}
\right.
$$

$(iv)$ From Theorem \ref{th2-11}, ${\rm la}_{k}(HS_4)\geq 4$. From
Observation \ref{obs1-1}, we have
$$
{\rm la}_{k}(HS_4)\leq {\rm la}_{k}(HP_4)+{\rm la}_{k}(HD_4)\leq
\left\{
\begin{array}{ll}
13 &\mbox {\rm if} \ k=1,\\
10 &\mbox {\rm if}\ k=2,\\
10 &\mbox {\rm if} \ k=3,\\
7 &\mbox {\rm if} \ k\geq 4.\\
\end{array}
\right.
$$
\end{pf}

\end{document}